\documentclass{amsart}

\usepackage{amssymb}
\usepackage[all]{xy}

\newtheorem{thm}{Theorem}%[section]
\newtheorem{mainthm}[thm]{Main Theorem}

\newtheorem{lem}[thm]{Lemma}
\newtheorem{cor}[thm]{Corollary}

\newtheorem{prop}[thm]{Proposition}

\theoremstyle{definition}
\newtheorem{defn}[thm]{Definition}

\newtheorem{say}[thm]{}
\newtheorem{exmp}[thm]{Example}

\newtheorem{ques}[thm]{Question}    %!!!!!!!!!!!!!!!!!!!!

\newtheorem{rem}[thm]{Remark}          

\newtheorem*{ack}{Acknowledgments}      % \renewcommand{\theack}{} 

\newtheorem{defn-thm}[thm]{Definition--Theorem}  %!!!!!!!!!!!!!!!!!!!!!!!!
\newtheorem{defn-lem}[thm]{Definition--Lemma}  %!!!!!!!!!!!!!!!!!!!!!!!!
\theoremstyle{remark}

\let \cedilla =\c
\renewcommand{\c}[0]{{\mathbb C}}  

\renewcommand{\o}[0]{{\mathcal O}} 
\newcommand{\z}[0]{{\mathbb Z}}

\renewcommand{\a}[0]{{\mathbb A}}

\newcommand{\p}[0]{{\mathbb P}}

\newcommand{\q}[0]{{\mathbb Q}}
\newcommand{\map}[0]{\dasharrow}
\newcommand{\qtq}[1]{\quad\mbox{#1}\quad}

\newcommand{\mult}[0]{\operatorname{mult}}

\newcommand{\supp}[0]{\operatorname{Supp}}    
    
\newcommand{\codim}[0]{\operatorname{codim}}

\newcommand{\sing}[0]{\operatorname{Sing}}    
\newcommand{\ex}[0]{\operatorname{Ex}}

\newcommand{\rup}[1]{\lceil{#1}\rceil}
\newcommand{\rdown}[1]{\lfloor{#1}\rfloor}

\newcommand{\onto}[0]{\twoheadrightarrow}

\newcommand{\simq}[0]{\sim_{\q}}

\newcommand{\tsum}[0]{\textstyle{\sum}}

\newcommand{\bir}[0]{\operatorname{Bir}}

\def\into{\DOTSB\lhook\joinrel\to}

\def\loccoh#1.#2.#3.#4.{H^{#1}_{#2}(#3,#4)}

\DeclareMathAlphabet{\mathchanc}{OT1}{pzc}%
                                {m}{it}

\usepackage[all]{xy}\xyoption{dvips}

\newcommand{\tprod}[0]{\textstyle{\prod}}

\begin{document}
\bibliographystyle{amsalpha}

%\today

\title[Rigidity theorem]{The rigidity theorem of\\
Fano--Segre--Iskovskikh--Manin--Pukhlikov--Corti--Cheltsov--de{\,}Fernex--Ein--Musta\cedilla{t}\u{a}--Zhuang}

  \author{J\'anos Koll\'ar}

\begin{abstract} We prove that $n$-dimensional smooth hypersurfaces of degree $n+1$ are superrigid. Starting with the work of Fano in 1915, 
the  proof of this Theorem  took 100 years and a dozen  researchers to construct. 
 Here I give complete proofs, aiming to use only basic knowledge of algebraic geometry and some Kodaira type vanishing theorems.

\end{abstract}

  \maketitle

\tableofcontents

The classification theory of algebraic varieties---developed by Enriques for surfaces and extended by Iitaka and then Mori to higher dimensions---says that every variety can be built from 3 basic types:
\begin{itemize}
\item (General type)  the canonical class $K_X$ is ample,
\item (Calabi-Yau)   $K_X$ is trivial and
\item (Fano)  $-K_X$ is ample.
\end{itemize}
Moreover, in the Fano case the truly basic ones are those that have
{\it class number} equal to 1. That is, every divisor $D$ on $X$ is
linearly equivalent to a (possibly rational) multiple of $-K_X$. 

If 2 varieties $X_1, X_2$  on the basic type   list are birationally equivalent then they have the same type. In the general type case they are even isomorphic and in 
 the Calabi-Yau case  the 
possible  birational maps   are reasonably well understood, especially for 3--folds, see   \cite{MR986434, k-etc}. 

By contrast, Fano varieties are sometimes birationally equivalent in quite unexpected ways and  the Noether--Fano method aims to  understand what happens. 

\begin{defn}[Weak rigidity and  superrigidity] I call a Fano variety $X$ with class number 1  {\it weakly rigid} if it is not birational to any other Fano variety $Y$ with class number 1, and
{\it weakly superrigid} if 
every birational map  $\Phi:X\map Y$ to another Fano variety $Y$ with class number 1  is an isomorphism.

 The adjective ``weakly'' is not standard; it allows us to define these notions without first discussing terminal singularities and  Mori fiber spaces.
The  definitions  of {\it  rigid} and {\it  superrigid} are similar, but allow  $Y$ to have terminal singularities and to be  a Mori fiber space; see \cite{MR1730270, MR2195677}.  

There are  many Fano varieties, especially in dimensions 2 and 3,
that are rigid but not superrigid.
% not birational to any other Fano variety but there are
% many birational self-maps  $\Phi:X\map X$ that are not isomorphisms.
% Such Fano varieties are called {\it rigid.}
Superrigidity is the more basic notion, though, in dimension 3, the theory of rigid Fano varieties is very rich.
\end{defn}

The aim of these notes is  to explain the proof of the following theorem. From now on we work over a field of characteristic 0. It is not important, but we may as well assume that it is algebraically closed.

\begin{mainthm}\label{superrigid.thm}
Every smooth hypersurface $X\subset   \p^{n+1}$ of dimension $n\geq 3$ and of degree $n+1$ is weakly superrigid.
\end{mainthm}

The proofs in the theory are designed to prove superrigidity, and
the optimal version of  Theorem~\ref{superrigid.thm} says that a smooth Fano hypersurface $X\subset   \p^{n+1}$  of dimension $\geq 3$ is  superrigid if and only if  $\deg X=n+1$; see \cite{MR2195677, puk-book}. The proof of this version needs only some new definitions and   minor changes in Step \ref{main.steps}.1.

% Equivalently, 
%  every birational map  $\Phi:X\map Y$ to another smooth Fano variety $Y$ with Picard number 1  is an isomorphism.

A smooth hypersurface $X\subset   \p^{n+1}$ of dimension $n\geq 3$ has class number 1 by 
Lefschetz's   theorem (see \cite{MR0033557} or \cite[p.156]{gri-har}).  
If $n=2$ then $X$ is a cubic surface, hence it has class number 7.
However, if the base field is not algebraically closed, it frequently  happens that $X$ has  class number 1, in which case it is weakly rigid but usually not weakly superrigid  by
\cite{MR0009471}; see  
 \cite[Chap.2]{ksc} for a modern treatment.

\begin{say}[The history of Theorem~\ref{superrigid.thm}] 

The first similar result is Max~Noether's description of all birational maps $\p^2\map \p^2$ \cite{MR1509694}. Noether's method formed the basis of all further developments. It was used by Segre to study birational maps of cubic surfaces over arbitrary fields \cite{MR0009471}, and  later generalized by
Manin and Iskovskikh to a birational theory of all
del~Pezzo surfaces and 2-dimensional conic bundles \cite{M66, isk79a}. 

 Theorem~\ref{superrigid.thm} was first stated by
Fano for 3-folds \cite{Fano-1908, Fano-1915}. His arguments contain many of the key ideas, but they also have  gaps. 
I call this approach the {\it   Noether--Fano method.}
The first complete proof  for 3-folds, along the lines indicated by Fano, is in
 Iskovskikh-Manin \cite{MR0291172}.  Iskovskikh and his school used this method to prove similar results for 
many other 3-folds, see \cite{MR537686,  Sar81, MR1668579, MR1859707}.
This approach  was gradually  extended to higher dimensions  by
Pukhlikov  \cite{MR870730, MR1650332, MR1970356} and  Cheltsov \cite{MR1805602}.
These results were complete up to dimension 8, but needed some additional general position assumptions in higher dimensions.
 A detailed  survey of this direction  is  in  \cite{puk-book}.

The theory of Fano varieties may be the oldest topic of higher dimensional birational geometry, but for a long time it grew almost independently of Mori's  Minimal Model Program. 
The Fano--Iskovskikh classification of Fano 3--folds  using extremal rays and flops was first  treated by Mori \cite{mori-83} and later improved by Takeuchi \cite{takeuchi}.

The  Noether--Fano method and the  Minimal Model Program  were brought  together   by 
Corti \cite{MR1311348}. Corti's technique  has been very successful in many cases, especially for 3--folds; see \cite{MR1798978} for a detailed study and \cite[Chap.5]{ksc} for an introduction. 
However, usually one needs some special tricks to make the last steps work, and a good higher dimensional version proved elusive for a long time.

New methods involving multiplier ideals were introduced by 
 de~Fernex-Ein-Musta\cedilla{t}\u{a} \cite{MR1981899}; these led to a more streamlined proof that worked up to dimension 12.   
The proof of Theorem~\ref{superrigid.thm} was finally completed by de~Fernex \cite{MR3455160}.

The recent paper of Zhuang \cite{zhuang}  makes the final step of  the Corti approach much easier in higher dimensions.
 The papers \cite{df17, MR3598821, sti-zhu, zhuang, liu-zhuang} contain more general results and  applications. 

The name 
Fano--Segre--Iskovskikh--Manin--Pukhlikov--Corti--Cheltsov--\-de{\,}Fernex--Ein--Musta\cedilla{t}\u{a}--Zhuang theorem was chosen to give credit to all those with a substantial contribution to the  proof, though this under emphasizes the major contributions of Fano, Iskovskikh and Pukhlikov.

The methods apply to many other Fano varieties for which $-K_X$ is a generator of the class group; see \cite{puk-book, zhuang} for  several examples. One of the big challenges is to understand what happens if 
 $-K_X$ is a multiple of the generator, see  \cite{MR3552255}.
\end{say}

\medskip

{\bf Open problems about  hypersurfaces}
\medskip

The following  questions are stated in the strongest forms that are consistent with the known examples. I have no reasons to believe that the answer to either of them is positive, and there may well be rather simple counter examples. As far as I know, there has been very little work on low degree hypersurfaces beyond cubics in dimension 4.

\begin{ques} Is every  smooth   hypersurface of degree $\geq 4$ non-rational? 
\end{ques}

Non-rationality of a smooth   hypersurface $X\subset \p^{n+1}$ is obvious if $\deg X\geq n+2$. For very general hypersurfaces of degree $\geq \tfrac23 n+3$, non-rationality was proved in \cite{Kollar95a}, a  major improvement by Schreieder \cite{schreieder} shows this for $\deg X\geq \log_2n +2$.

\begin{ques} \label{weak.neq.sup.hyp.ques} 
 Is every  smooth   hypersurface of degree $\geq 5$ weakly superrigid?
\end{ques}

Here $\geq 5$ is necessary since 
there are  some smooth quartics with nontrivial birational maps.

\begin{exmp}\label{quart.2n.2L.inv.prop}
 Let $X\subset \p^{2n+1}$ be a quartic hypersurface that contains 2 disjoint linear subspaces $L_1, L_2$ of dimension $n$. For every $p\in \p^{2n+1}\setminus (L_1\cup L_2)$ there is a unique line $\ell_p$ through $p$ that meets both $L_1, L_2$. 
This line meets $X$ in 4 points, two of these are on $L_1, L_2$. If $p\in X$ then this leaves a unique 4th intersection point, call it $\Phi(p)$. 
Clearly $\Phi$ is an involution which is not defined at $p$ if either
$p\in L_1\cup L_2$ or if $\ell_p\subset X$. 
\end{exmp}

\begin{ack} These notes are based on lectures given at the AGNES Conference at Rutgers and at Northwestern University in April and May, 2018. The hospitality of Northwestern University, and especially of M.~Popa,  gave an ideal time to write them up in expanded form. A received many helpful comments and references from N.~Alon, 
 V.~Cheltsov, A.~Corti, L.~Ein  and Y.~Liu.
Several e-mails from  T.~de~Fernex helped me to  clarify some of the key points and to improve the paper considerably.
 
Partial  financial support    was provided  by  the NSF under grant number
 DMS-1362960 and by the Nemmers Prize of Northwestern University.
\end{ack}

\section{Rigidity and superrigidity, an overview}

In the following outlines and in  subsequent Sections 
I aim to   put the pieces together and write down a simple proof of the superrigidity of  smooth hypersurfaces
$V^n_{n+1}\subset \p^{n+1}$,  where most steps are either easy or are direct applications of some general principle of the  Minimal Model Program. 

The key  notion we need is  canonical and log canonical
pairs involving linear systems. 

\begin{defn}[Log resolution]\label{log.res.say} Assume that we have a  variety $X$, a (not necessarily complete) linear system $|M|$ on  $X$ and a divisor $D$ on $X$. 
A {\it log resolution} of these data is
a proper birational morphism $\pi:X'\to X$ 
such that 
\begin{enumerate}
\item $X'$ is smooth,   
\item $\pi^*|M|=|M'|+B$ where $|M'|$ is base-point free and $B$ is the fixed part of $\pi^*|M|$, and 
\item $B+\pi^{-1}_*(D)+\ex(\pi)$ is a simple normal crossing divisor.
\end{enumerate}
(Here $\ex(\pi)$ denotes the exceptional set of $\pi$,
$\pi^{-1}_*(D) $ denotes the birational transform of $D$ and {\it simple normal crossing} means that the irreducible components are smooth and they intersect transversally. The adjective ``log'' loosely refers to condition (3).)

The existence of log resolutions was proved by Hironaka; see
\cite[Chap.3]{k-res} for a recent treatment.
\end{defn}

\begin{defn}[Canonical, log canonical, etc.]\label{can.logcan.defn} Let $X$ be a smooth variety and $|M|$ a linear system on  $X$. Let $\pi:X'\to X$ be a log resolution of $|M|$ 
as in Definition~\ref{log.res.say}. 
  Write  $\pi^*|M|=|M'|+B$ where $|M'|$ is base-point free and $B$ is the fixed part of $\pi^*|M|$, and $K_{X'}\sim \pi^*K_X+E$ where $E$ is effective and $\pi$-exceptional. For any nonnegative rational number $c$ we can thus formally write
$$
K_{X'}+c|M'|\simq \pi^*\bigl(K_X +c|M|\bigr)+(E-cB),
\eqno{(\ref{can.logcan.defn}.1)}
$$
where $A_1\simq A_2$ means that $N\cdot A_1$ is linearly equivalent to $N\cdot A_2$ for some $N>0$. 

A pair  $\bigl(X, c|M|\bigr)$ is called  {\it canonical}  (resp.\ 
  {\it log canonical}) if every divisor appears in
$E-cB$ with coefficient $\geq 0$  (resp.\ $\geq -1$).
This is independent of the log resolution \cite[2.32]{km-book}. 

Note that if $r>0$ is an integer then $\bigl(X, c|M|\bigr)$ is canonical  (resp.\  log canonical) iff $\bigl(X, \tfrac{c}{r}|rM|\bigr)$ is.
(Keep in mind that  $|rM|$ is the linear system spanned by sums of the form
$M_1+\cdots+M_r$ where $M_i\in |M|$.)
Thus we can always restrict to dealing with pairs  $\bigl(X, c|M|\bigr)$  where $c< 1$; this is frequently convenient.

If $|M|$ is base point free then $B=0$, thus
 $\bigl(X, c|M|\bigr)$ is  canonical for any $c$. In all other cases,
 $\bigl(X, c|M|\bigr)$ is  canonical (resp.\  log canonical) for small values of $c$ but not for large values. (The transitional value of $c$ is called the
canonical (resp.\  log canonical)  {\it threshold.})
Roughly speaking, small threshold corresponds to very singular base locus.

Write  $E-cB=\sum_i a_i E_i$. Then $E_i$ is  called a  {\it non-canonical divisor}
(resp.\  a  {\it non-log-canonical divisor}) of   $\bigl(X, c|M|\bigr)$ iff  $a_i< 0$  (resp.\ $a_i< -1$).  
The  corresponding image 
$\pi(E_i)\subset X$   is a   {\it  non-canonical center}
(resp.\  {\it non-log-canonical center})   of $\bigl(X, c|M|\bigr)$. 
These centers are always contained in the base locus of $|M|$. 
It is not very important for us, but, as we run through all log resolutions and all divisors on them,
we might get infinitely many non-(log)-canonical centers, however their union is the  closed subset
$\textstyle{\bigcup_{a_i<0}}\pi(E_i)$ (resp.\
$\textstyle{\bigcup_{a_i<-1}}\pi(E_i)$);
 see  \cite[2.31]{km-book}.

Using Remark~\ref{div.ls.id.rem}, these formulas also define the above notions for pairs $(X, \Delta)$ where $\Delta$ is an effective  divisor and 
pairs $\bigl(X, I^c\bigl)$ where $I$ is an ideal sheaf.
\end{defn}

\begin{rem}[Divisors, linear systems and ideal sheaves]
\label{div.ls.id.rem}
Much of the  Minimal Model Program literature works with pairs
$(X,\Delta)$ where $\Delta$ is a divisor (with rational or real coefficients),
see \cite{km-book, kk-singbook}. For rigidity questions, the natural object seems to be a pair $\bigl(X, c|M|\bigl)$ where  $|M|$ is a linear system and $c$ is a rational or real coefficient.  It is easy to see that if
$c\in [0,1)$  (which will always be the case for us) and  $D\in |M|$ is a general divisor then the definitions and theorems for $\bigl(X, c|M|\bigl)$ and $\bigl(X, cD\bigl)$ are equivalent. 

As we noted in Definition~\ref{can.logcan.defn}, working with $\bigl(X, |M|\bigr)$ is  equivalent to working with $\bigl(X, \tfrac12 |2M|\bigr)$, but one version may give a clearer picture than the other. 
  As an illustration, consider the linear system $|\lambda x+\mu( x-y^r)=0|$ in the plane. A general member of it is a smooth curve and the role of the $y^r$ term is not immediately visible.
By contrast the linear system $|2M|$ is  $|\lambda x^2+\mu x(x-y^r)+\nu (x-y^r)^2=0|$, its general member is (after a local analytic  coordinate change) of the form $(x^2-y^{2r}=0)$. 
Now we see both the original smoothness  (since $x^2$ is there) and the order of tangency between two members of $|M|$  (shown by $y^{2r}$). 
While computationally this is not important, conceptually it seems clearer  that  information about intersections of 2 divisors in $|M|$ is now visible on individual divisors in $|2M|$.

Let $X$ be an affine variety and $I\subset \o_X$ an ideal sheaf.
Many authors, for example
\cite{MR2047683, laz-book}, work with  pairs $\bigl(X, I^c\bigl)$ where  $c$ is viewed as a formal exponent. 
% (If $c$ is a natural number, one can identify the formal $I^c$ with the $c$th power of $I$.) 
If $I$ is
 generated by global sections $g_1,\dots, g_m$, we can consider the linear system  $|M|:=|\sum \lambda_ig_i=0|$. Again we find that  the definitions and theorems for $\bigl(X, c|M|\bigl)$ and $\bigl(X, I^c\bigl)$ are equivalent. 

Here I  follow the language of linear systems, since this seems best suited to  our current aims. I will also always assume that $c$ is rational. This is always the case in our applications and makes some statements  simpler. However, it does not cause any essential difference at the end.
\end{rem}

We discuss the canonical and  log canonical property of linear systems in detail in Section~\ref{sec.12}. For now we mainly need to know that
canonical means  mild singularities and  log canonical means somewhat worse singularities. In some sense the main question  of the theory was how to describe these properties  in terms of other, better understood, measures of singularities.

\begin{say}[Main steps of the proof]\label{main.steps} The proof can be organized into 6 fairly independent steps. Roughly speaking, Steps 1 and  2  are essentially in the works of Fano, at least for 3--folds. Steps 3 and 4 are substantial reinterpretations of the classical ideas while Steps 5 and  6 give a  new way of finishing the proof. 
\medskip

{\bf Notation.} For the rest of this section  I write $Y$ for a smooth, projective variety, 
$X$ for a smooth, projective Fano variety with class number 1 and $V$ (or $V_{n+1}$ or $V^n_{n+1}$) for a smooth hypersurface
$V^n_{n+1}\subset \p^{n+1}$ of degree $n+1$ and of dimension $n\geq 3$. The base field has characteristic 0.

\medskip
{\bf Step \ref{main.steps}.1} (Noether-Fano criterion, Section~\ref{sec.10})  {\it A smooth Fano variety $X$ of class number 1  is (weakly) superrigid if
for every movable linear system  $|M|\subset |-mK_X|$ the pair
 $\bigl(X, \tfrac{1}{m} |M|\bigr)$ is canonical.}

\medskip
{\it Comments  \ref{main.steps}.1.1.} Movable means that there are no fixed components, some authors use {\it mobile} instead. If $\dim X=2$ then $\bigl(X, \tfrac{1}{m} |M|\bigr)$ is not canonical iff $\mult_x|M|> m$ for some point  $x\in X$ by Lemma~\ref{main.steps.3.a.pf}; this equivalence made Noether's and Segre's  proofs work well. If $\dim X=3$ then Fano tried to prove that if  $\bigl(X, \tfrac{1}{m} |M|\bigr)$ is not canonical  then either $\mult_C|M|> m$ for some curve $C\subset X$ or 
$\mult_x|M|> 2m$ for some point $x\in X$. Fano understood that the latter  condition for points is not right, one needs instead only a consequence of it:
The local intersection number at $x$ is $(M\cdot M\cdot H)_x>4m^2$, where $H$ is a hyperplane through $x$.  
In higher dimensions it does not seem possible to define canonical in terms of just multiplicities and intersection numbers, this is one reason why the above form of Step \ref{main.steps}.1  was  established only in  \cite{MR1311348}.
We   prove  Step \ref{main.steps}.1 in Theorem~\ref{noe.fano.ineq.thm}. 
Although historically the notion of  ``canonical'' was first defined starting from varieties of general type  (see \cite{r-c3f, km-book}),  the Noether-Fano criterion leads to the exact same notion.

If $X$ is not (weakly) superrigid then there is a  movable linear system  $|M|\subset |-mK_X|$ such that
 $\bigl(X, \tfrac{1}{m} |M|\bigr)$ is not canonical, thus it has some
non-canonical divisors and centers as in Definition~\ref{can.logcan.defn}.
(The  ``worst''  non-canonical centers  are  called {\it maximal centers}  by the Iskovskikh school.)
From now on  we  focus  entirely on understanding movable linear systems and their possible non-canonical centers on $X$. 
There are 2  persistent problems that we encounter.
\begin{itemize}
\item We can usually bound the {\it multiplicities} of $|M|$, but there is a gap---growing with the dimension---between multiplicity and the canonical property.
\item We are  better at understanding when a pair is 
{\em log} canonical, instead of canonical.
\end{itemize}

While we try to make statements about arbitrary Fano varieties, at some point we need to use special properties of the $V^n_{n+1}$.
The following bounds, going back to Fano and Segre, were put into final form by Pukhlikov \cite[Prop.5]{MR1970356}   and later generalized by Cheltsov \cite[Lem.13]{MR2171057} and Suzuki \cite[2.1]{MR3598821} to complete intersections. 
%See \cite[5.36]{ksc} for the simplest version.

\medskip
{\bf Step \ref{main.steps}.2} (Multiplicity bounds, Fano, Segre, Pukhlikov, Section~\ref{sec.11})   {\it   Let  $Y\subset \p^{n+1}$ be a smooth hypersurface and $|H|$ the hyperplane class on $Y$.
Let $D\in |mH|$ be a divisor, $|M|\subset |mH|$  a movable linear system and 
$Z\subset Y$  an irreducible subvariety.  
\begin{enumerate}
\item[(a)] If $\dim Z\geq 1$ then $\mult_{Z}D\leq m$.
\item[(b)] If $\dim Z\geq 2$ then $\mult_{Z}(M\cdot M)\leq m^2$,
\end{enumerate}
where $M\cdot M$ denotes the intersection of 2 general members of $|M|$.
}

\medskip
{\it Comments  \ref{main.steps}.2.1.} 
Note that (\ref{main.steps}.1) works with
$|M|\subset |-mK_X|$  and (\ref{main.steps}.2) with $|M|\subset |mH|$.
The two match up iff $-K_X\sim H$; the latter holds for $X=V^n_{n+1}$, the case that  we are  considering. In general, the method works best for those Fano varieties where every divisor is an integral  multiple of $-K_X$  (up to linear equivalence). 
\medskip

Next we need to understand the relationship between the 
multiplicity bounds in Step~\ref{main.steps}.2 and the  canonical property. This is rather 
 easy  for  (\ref{main.steps}.2.a). Combining it with  (\ref{main.steps}.3.a)  we get that   $\bigl(Y, \frac1{m}D\bigr)$ is canonical, except at a finite point set $P\subset Y$.
We already mentioned this in Comments~\ref{main.steps}.1.1;
 see Lemma~\ref{main.steps.3.a.pf} or  \cite[6.18]{ksc} for proofs.  

Relating  (\ref{main.steps}.2.b) to the canonical  property was less obvious; it was done by Corti \cite[3.1]{MR1798984}   (see also  \cite[Sec.6.6]{ksc}),
 then very much generalized by de~Fernex-Ein-Musta\cedilla{t}\u{a}  \cite{MR2047683}  and  sharpened by Y.~Liu \cite{liu16}.
%The starting point is a connection with the  log canonical  property.

\medskip
{\bf Step~\ref{main.steps}.3} (Non-(log)-canonical points and multiplicity, Corti, Sections~\ref{sec.12}--\ref{sec.13})   {\it Let  $|M|$ be a movable linear system on a smooth variety $Y$.
\begin{enumerate}
\item[(a)]  
If $\bigl(Y, \frac1m |M|\bigr)$ is not canonical at ${\mathbf p}\in Y$ then $\mult_{\mathbf p}|M|> m$.
\item[(b)]  If $\bigl(Y, \frac1m |M|\bigr)$ is not log canonical at ${\mathbf p}\in Y$ then $\mult_{\mathbf p}(M\cdot M)> 4m^2$.
\end{enumerate}
}

\medskip
{\it Comments  \ref{main.steps}.3.1.}  Both of these bounds are sharp as shown by the examples   $$
\bigl(\a^2, \tfrac1m |\lambda x^{m+1}+\mu y^{m+1}=0|\bigr)\qtq{and}
\bigl(\a^2, \tfrac1m |\lambda x^{2m+1}+\mu y^{2m+1}=0|\bigr),
$$ 
which have a 
 non-canonical (resp.\ non-log-canonical) center at the origin. Surprisingly, part (a) can not  be improved for non-log-canonical centers, as shown by 
$$
\bigl(\a^2, \tfrac1m |\lambda x^{m+1}+\mu y^{(m+1)^2}=0|\bigr).
$$ (This can be computed by hand or  see (\ref{lcthreshold.curve.formula}).)

Using  a---by now standard---method  called {\it inversion of adjunction,}  which we discuss in Section~\ref{sec.13}, 
both parts   follow from   claims about linear systems on algebraic surfaces:
\medskip

{\it Claim \ref{main.steps}.3.2.}  Let  $|M|$ be a movable linear system on a smooth surface $S$.
\begin{enumerate}
\item[(a)] If $\bigl(S, c |M|\bigr)$ is not  canonical at $s\in S$ then $\mult_{s}|M|> 1/c$.
\item[(b)]  If $\bigl(S, c |M|\bigr)$ is not log canonical at $s\in S$ then $\mult_{s}(M\cdot M)> 4/c^2$.
\end{enumerate}

% There are 2 ways to use (\ref{main.steps}.3.b). The first uses that we get  a lower bound on some intersection numbers.  As an  illustration, we consider  quartic threefolds.

\medskip

It would be very nice to continue the claims (\ref{main.steps}.3.a--b) to stronger and stronger inequalities for higher codimension non-log-canonical centers. This was done in \cite{MR2047683}. This is very useful if by chance the base locus of $|M|$ has codimension $>2$. However, in many cases the 
base locus of $|M|$ has codimension $2$ and it is not easy to apply the estimates of \cite{MR2047683} directly.

% We need a new way of exploiting the tension between the estimates
% (\ref{main.steps}.2.b) and (\ref{main.steps}.3.b). The next step gives  the other way of  using (\ref{main.steps}.3.2.b).

\medskip

Fano always aimed to reduce questions about Fano 3-folds 
(for him these meant $X\subset \p^n$ such that $-K_X\sim H$) to their hyperplane sections. These  are K3 surfaces, whose geometry was quite well understood. In higher dimensions, the  hyperplane sections are Calabi-Yau varieties, whose geometry is much less known.  Thus the modern focus is on the change of the singularities as we  restrict a linear systems to a hyperplane section.

\medskip
{\bf Step \ref{main.steps}.4} (Cutting by hyperplanes, Section~\ref{sec.13}) {\it  Let  $|M|$ be a movable linear system on a smooth variety $Y$.  Fix a point $p\in Y$ and let $W$  be a general member of a very ample linear system  $|H|$ that passes through $p$.
Then 
\begin{enumerate}
\item[(a)]  If  $\bigl(Y, c |M|\bigr)$ is (log) canonical outside a closed subset $Z\subset Y$ then 
$\bigl(W, c |M|_W\bigr)$ is (log) canonical outside  $\{p\}\cup (Z\cap W)$. 
\item[(b)] If  $\bigl(Y, c |M|\bigr)$ is not log canonical at $p$ then
$\bigl(W, c |M|_W\bigr)$ is also not log canonical at $p$.
\item[(c)] If $p$ is a non-canonical center of  $\bigl(Y, c |M|\bigr)$ then
$p$ is a non-log-canonical point of $\bigl(W, c |M|_W\bigr)$.
\end{enumerate}
}
\medskip

{\it Warning.} Note that in (c) the point $p$ needs to be a 
 non-canonical {\em center} on $Y$ 
(see Definition~\ref{can.logcan.defn}) and then it is 
 a non-{\em log}-canonical point on $W$. 

\medskip
{\it Comments  \ref{main.steps}.4.1.}  
The multiplicity versions of these go back to Bertini and Fano, 
  but the above form of (c)   may have been first made explicit in \cite{MR1798984}.  By now these are special cases of the theory of {\it adjunction} for log canonical pairs, we discuss this in  Section~\ref{sec.13}.

Note that 
cutting by a hyperplane has a very curious effect on the singularities. 

If $\bigl(X, c|M|\bigr)$ is a canonical (resp.\ log canonical) pair then its restriction  to a general member of a  base point free linear system is still  canonical (resp.\ log canonical); this is an easy Bertini-type theorem, see (\ref{bertini.say}.1).  Applying this to  $X:=Y\setminus Z$ gives (a). 
Part (b) is quite a bit harder to prove but it fits  the general pattern that singularities do not get better by cutting with a hyperplane. 

The surprising part is (c) which says that the singularity is  made   {\em worse} by restriction to a general hypersurface through a non-canonical center.  This is in marked contrast with multiplicity, which is preserved by such restrictions. We discuss this in Section~\ref{sec.13}.

% However, if  $p\in X$ is a  non-canonical center of  $\bigl(X, c|M|\bigr)$ and
% $W$ passes through $p$ then  $\bigl(W, c|M|_W\bigr)$ is not {\em log}  canonical at $p$. That is, 

The   first application of these ideas is the following.
\medskip

{\it  Rigidity of quartic threefolds, Corti's variant. \ref{main.steps}.4.2.} Let 
 $X\subset \p^4$ be a  smooth quartic 3-fold and assume that 
we have a linear system $|M|\subset |mH|$ such that 
 $\bigl(X, \frac1m |M|\bigr)$ is not canonical.  
One dimensional non-canonical centers are excluded by (\ref{main.steps}.2.a) and (\ref{main.steps}.3.a).
If $x\in X$  is a 0-dimensional non-canonical center  then
let $W$ be a general hyperplane section passing through $x$.
Then $x$ is a non-log-canonical center of $\bigl(W, \frac1m |M|_W\bigr)$
by (\ref{main.steps}.4.c), hence  the local intersection number 
 $(M\cdot M \cdot W)_x$ is $>4m^2$  by (\ref{main.steps}.3.b). Therefore
$(M\cdot M\cdot W)>4m^2$.  On the other hand,  $M\sim mH$ and hence $(M\cdot M\cdot W)=4m^2$, a contradiction.
Thus smooth  quartic 3-folds are weakly superrigid. \qed
\medskip

More generally,   the method described so far works well  if $-K_X$ generates the class group and  $(-K_X)^n\leq 4$. 
Among hypersurfaces in $\p^{n+1}$, this holds only for the quartic 3-folds. However there are smooth hypersurfaces in weighted projective spaces with these properties. For example, fix $r>1$ and  let $X$ be a smooth hypersurface of degree $4r + 2$ and dimension $2r$ in the weighted projective space $\p(1^{2r}, 2, 2r + 1)$ (the notation means that we have $2r$ coordinates of weight 1, see \cite[3.48]{ksc} for an introduction). 
Then  $-K_X\sim H$ and $(-K_X)^{2r}=1$. With small changes the method  proves that they are superrigid; see \cite[5.22]{ksc} for details. \qed
\medskip
 
% From now on we work with $\bigl(W, \tfrac1m |M|_W\bigr)$ and try to reach a contradiction.

In the above proof we have used (\ref{main.steps}.3.b) on $W\in |H|$.
The new idea of \cite{MR3455160} is to use it directly on $Y$.

\medskip
{\bf Step \ref{main.steps}.5} (Doubling the linear system, de~Fernex,  Sections~\ref{sec.12}--\ref{sec.13}) {\it Instead of working only with 
 $\bigl(Y, \tfrac1m |M|\bigr)$, we should focus on the interaction between 
 $\bigl(Y, \tfrac1m |M|\bigr)$ and $\bigl(Y, \tfrac1m |2M|\bigr)$.}

In order to contrast the 2 cases,  let  $Y\subset \p^{n+1}$ be a smooth hypersurface and  $|M|\subset |mH|$  a movable linear system.
Combining (\ref{main.steps}.2.a) with (\ref{main.steps}.3.a) and  (\ref{main.steps}.2.b) with (\ref{main.steps}.3.b) gives the following.
 {\it 
 \begin{enumerate}
\item[(a)] $\bigl(Y, \tfrac1m |M|\bigr)$ is   canonical outside a finite set of points  $P\subset Y$ and 
\item[(b)] $\bigl(Y, \tfrac1m |2M|\bigr)$ is  log canonical outside a finite set of curves  $C\subset Y$.
\end{enumerate}
}

%  Note the key shift of having
% $|2M|$ instead of $|M|$ in part (b). 

One should think of these  as saying that  $|M|$ and 
 $|2M|$ are very singular at $p$ but less  singular almost everywhere else. 
A key insight of \cite{MR3455160} is that (b) is much stronger than (a). 
In order to understand this, let us see how one can use the information provided by (\ref{main.steps}.5.a--b). 
 
% This is quite important in applications. Usually we can get a lot of information if we have an isolated ``very singular'' point $p$ but much less if $p$  is a limit of other singular points that are  ``almost as singular.''

\medskip
{\it Comments  \ref{main.steps}.5.1.} Let $Y$ be a smooth, projective variety, $H$ an ample divisor and $\Delta\simq H$ a 
 $\q$-divisors with an isolated non-log-canonical center at a point $p\in Y$. The observation that this leads to a global section of $\o_Y(K_Y+H)$ that does not vanish at $p$ has been  an important  ingredient  of the  Kawamata--Reid--Shokurov approach to the cone theorem (cf.\ \cite[Chap.3]{km-book})  and is  central in the works around Fujita's conjecture  (cf.\ \cite[Secs.5--6]{k-pairs} or \cite[Sec.10.4]{laz-book}). In all these applications the aim is to get at least 1 section that does not vanish at a given point. Although it was known that the process can be used to get several sections, this has not been the focus in the past.

At first sight, (\ref{main.steps}.5.a) is better suited to use this method.
If $\bigl(Y, \tfrac1m |M|\bigr)$ is   not canonical at some $p\in P$, then,
by (\ref{main.steps}.4.c), after restricting to a general hyperplane section $p\in W\subset Y$, we get  $\bigl(W, \tfrac1m |M|_W\bigr)$ that is   not  log canonical at  $p$ but is canonical outside $P\cap W$.
This leads to a section of $\o_W(K_W+H)$ that does not vanish at $p$. However, in our cases $\o_W(K_W+H)$ is very ample, so there is no  contradiction.

The problem seems to be that while we have been thinking of canonical as ``much better'' than log canonical, from the numerical point of view the difference seems small. We saw an instance of this in (\ref{main.steps}.3.a), where both the non-canonical and non-log-canonical cases yield the same  inequality; see   (\ref{main.steps}.3.1).

In  (\ref{main.steps}.5.b) we ``gain'' since $|M|$ is replaced by $|2M|$ but also
``lose'' since canonical is replaced by log canonical and the finite set of points $P$ is replaced by a finite set of curves  $C$. However, when we switch to a  hyperplane section $W\subset Y$, we focus on the 
  non-log-canonical property anyhow, and $C\cap W$ becomes a finite set of points.
So the ``losses'' do not matter at the end but
the ``gain''  stays with us. Thus we get that
\begin{enumerate}
\item[(a)]   $\bigl(W, \tfrac1m |M|_W\bigr)$  is   not  log canonical at  some $p\in P\cap W$, but even
\item[(b)] $\bigl(W, \tfrac1m |2M|_W\bigr)$  is    log canonical outside  $P\cap W$. 
\end{enumerate}
A fundamental claim of \cite{MR3455160} is that this in itself leads to a contradiction. While the argument at the end of \cite{MR3455160} is worded differently, unraveling the proofs of \cite[Lems.3--4]{MR3455160} 
gives a quadratic lower bound for  $h^0\bigl(W, \o_W(K_W+2H)\bigr)$, almost enough to get a contradiction without further work.
 Building on 
\cite{MR3455160}, a  key  observation of \cite{zhuang} is that a suitable modification of this method leads to an exponential lower bound and a quick numerical contradiction.

\medskip
{\it Remark  \ref{main.steps}.5.2.} Once the technical details are settled, we see that there is  lot of room in (\ref{main.steps}.5.b).  Namely, if we know only that, for some fixed $\epsilon>0$ and $d$, $\bigl(Y, \tfrac1m |(1+\epsilon)M|\bigr)$ is  log canonical outside a  subset of dimension $\leq d$, that is still enough to prove Theorem~\ref{superrigid.thm} for $n$ sufficiently large (depending on $\epsilon$ and $d$).

% We still have a problem with the lack of control along the curves $C$ in
% (\ref{main.steps}.5.b). We go around it  by taking a 
% hyperplane section.

 % \cite{zhuang} shows that  focusing on  the  multiplier ideal of a  multiplier ideal gives an exponential  lower bound.   

\medskip
{\bf Step \ref{main.steps}.6} (Zhuang, Section~\ref{sec.14})  {\it Let $Y$ be a smooth projective variety of dimension $d$
 and $L$  an ample  divisor on $Y$. Further let  $|M|\subset |mL|$ be a 
movable linear system and   $P\subset Y$ a finite (nonempty) subset of  $Y$. Assume that  
\begin{enumerate}
\item[(a)] $\bigl(Y, \tfrac1{m} |M|\bigr)$ is not log canonical at  some  $p\in P$, but
\item[(b)] $\bigl(Y, \tfrac1{m}|2M|\bigr)$ 
 is log canonical outside $P$. 
\end{enumerate} 
Then 
$$
h^0\bigl(Y, \o_Y(K_Y+2L)\bigr)\geq\tfrac12 3^d.
$$
}

\medskip
{\it Comments  \ref{main.steps}.6.1.}  One should think of this as saying that if $|M|$ is much more singular at a finite set of points than elsewhere then
the linear system $|K_Y+2L|$ is very large. I stated the case where we compare
the singularities of $ \tfrac1{m}|M|$ and $ \tfrac1{m}|2M|$, the complete version in \cite{zhuang}  also applies 
if we work with  $ c|M|$ and $ (c+\epsilon)|M|$ for some $\epsilon>0$. 

It is quite remarkable that there is also a rather easy converse.

Let $|L|$ be any linear system on $Y$ and $y\in Y$ a point.
If $\dim |2L|\geq \binom{3d}{d}$ then there is a linear subsystem  $|N|\subset |2L|$ that has multiplicity $>2d$ at $y$. In particular, $\bigl(Y, \frac1{2} |N|\bigr)$ is not
log canonical at $y$. 
As $d\to\infty$, $\binom{3d}{d}$ grows like $6.75^d$. 

Thus if $\dim |2L|\geq 6.75^d$ then we can find a linear system
$|N|$ that satisfies (\ref{main.steps}.6.a),   and usually also 
(\ref{main.steps}.6.b).  
Informally we can restate (\ref{main.steps}.6) as
\medskip

{\it Principle  \ref{main.steps}.6.2.
 There are  no accidental isolated singularities.}
\end{say}

\begin{say}[Proof of Theorem~\ref{superrigid.thm} using Steps \ref{main.steps}.1--6]\label{sup.rig.thm.pf}
Let $V\subset \p^{n+1}$ be a smooth hypersurface of degree $n+1\geq 4$. 
If $V$ is not weakly superrigid, then, 
 by Step \ref{main.steps}.1 we get  $\bigl(V, \tfrac1m |M|\bigr)$ that is not canonical. Thus Steps \ref{main.steps}.3--5 give a
$W=W^{n-1}_{n+1}\subset \p^{n}$ and  $|M|_W$ such that 
\begin{enumerate}%\setcounter{enumi}{2}
\item[(a)] $\bigl(W, \tfrac1{m}|M|_W\bigr)$ is not  log canonical at finitely many  points $P\subset W$, but 
\item[(b)] $\bigl(W, \tfrac1{m}|2M|_W\bigr)$ is  log canonical outside   $P$.
\end{enumerate}
By   Step \ref{main.steps}.6 this implies that
$$
h^0\bigl(W, \o_W(K_W+2H)\bigr)\geq \tfrac12 3^{n-1}.
$$
On the other hand 
$h^0\bigl(W, \o_W(K_W+2H)\bigr)=h^0\bigl(W, \o_W(2H)\bigr)=
h^0\bigl( \p^{n}, \o_{ \p^{n}}(2)\bigr)=\tbinom{n+2}{2}$, so
$$
\tbinom{n+2}{2}\geq\tfrac12 3^{n-1}.
%\eqno{(\ref{sup.rig.thm.pf}.4)}
$$
The left hand side is quadratic in $n$, the right hand side is exponential,
so for $n\gg 1$ this can not hold. (In fact, we have a lot of room, leading to many other  cases where the method applies in large dimensions; see \cite{zhuang}.)

By direct computation, we get a contradiction for $n\geq 5$, hence  we get the superrigidity of $V^n_{n+1}\subset \p^{n+1}$ for $n\geq 5$. 

One can improve the lower bound in (\ref{main.steps}.6) to $\tfrac12 3^d+\tfrac32$, and then for $n=4$ we get an equality  $\tbinom{6}{2}=15=\tfrac12 3^3+\tfrac32$. 
So there is no contradiction, but it is quite likely that a small change can make the proof work.
However, the   $n=3$ case does not seem to follow, but this was already treated in (\ref{main.steps}.4.2). \qed
\end{say}

\begin{say}[Attribution of the Steps] In rereading many of the contributions to the proof I was really struck by how gradual the progress was and how difficult it is to attribute various ideas to a particular author or paper. 

Fano's papers are quite hard to read, and some people who spent years on trying to learn from them came away with feeling that Fano got most parts of the proof wrong. Others who looked at Fano's works   feel that he had all the essential points right. 
In particular, the attribution of Step 2 has been controversial. 

I think of Corti's work \cite{MR1798984} as a major conceptual step forward, but some authors felt that it did not add anything new, at least initially. The idea of doubling the linear system is in retrospect already in \cite{MR1798984}, 
but the new viewpoint of de~Fernex \cite{MR3455160} turned out to be very powerful and, as we discussed in (\ref{main.steps}.5.1), 
 the latter contains many of the ingredients of Step~6.
 I had a hard time formulating  Steps~3--6 in a way that shows the differences between them  meaningfully while highlighting the new ideas of the main contributor.
 Nonetheless, at least in hindsight, each of the Steps represents a major new idea, though this was not always immediately understood.

No doubt several people will feel that my presentation is flawed in many ways.
Luckily the reader can consult the excellent survey
\cite{MR2195677} and   books
\cite{MR1798978, puk-book} for different viewpoints. 
\end{say}

\begin{say}[What is missing?] My aim was to write down a proof of
Theorem~\ref{superrigid.thm} that is short and focuses on the key ideas.
My preference is for steps that follow from general results and techniques of the MMP. 
Thus several important developments have been left out.

After proving rigidity for quartic 3--folds, the Russian  school went on to study other Fano 3--folds. They found that they are frequently 
rigid but not superrigid and the main question is how to find generators for $\bir(X)$.
 The contributions of Iskovskikh, Sarkisov, Pukhlikov and Cheltsov are especially significant.
These results and their higher dimensional extensions are surveyed in 
\cite{MR2195677, puk-book}.

The first major applications of the Corti method were also in dimension 3, see  \cite{MR1798978} for a survey and \cite{MR3019454} for a higher dimensional extension.

In our proof we need to understand 0-dimensional log canonical centers, but  the theory of arbitrary log canonical centers has been quite important in higher dimensional geometry.
The first  structure theorems were proved by Ambro \cite{ambro};
 see \cite[Chaps.4--5]{kk-singbook} and \cite{fuj-book} for later treatments and generalizations. 
%The theory was further developed in \cite{MR2047683}.
\end{say}

\section{The Noether-Fano method}\label{sec.10}

We start the proof of Theorem~\ref{superrigid.thm} by establishing 
Step \ref{main.steps}.1.

\begin{thm}[Noether-Fano inequality]\label{noe.fano.ineq.thm}
 Let  $\Phi:X\map X'$ be a birational map between  smooth Fano varieties of class number 1. Then
\begin{enumerate}
\item either $\Phi$ is an isomorphism,
\item or there is a movable linear system $|M|\subset  |-mK_X|$ for some $m>0$ on $X$ such that
$\bigl(X, \tfrac1m |M|\bigr)$ is not canonical. 
\end{enumerate}
\end{thm}

Proof. 
Let $Z$  be the normalization of the closure of the graph of $\Phi$ with projections
$p:Z\to X$ and $q:Z\to X'$.
Pick any base-point-free linear system $ |M'|\subset |-m'K_{X'}|$ and let
  $|M|:=\Phi^{-1}_*|M'|$ denote its birational transform on $X$. Set  $|M_Z|=q^*|M'|$.
Since the class number of  $X$ is 1,   $|M|\simq  -mK_X$ for some
$m>0$.  (If $m$ is not an integer, we replace $|M'|$ by a suitable multiple. Thus we may as well assume that $|M|\subset  |-mK_X|$.)  We define a $q$-exceptional divisor $E_q$ and  $p$-exceptional divisors $E_p, F_p$ by the formulas
$$                                                                               \begin{array}{ll}
K_Z=q^*K_{X'}+E_q, & |M_Z|=q^*|M'|\qtq{and}\\                                          
K_Z=p^*K_X+E_p,  & |M_Z|=p^*|M|-F_p.
\end{array}                                           
\eqno{(\ref{noe.fano.ineq.thm}.3)}
$$
Since $X', X$ are smooth,  $E_q, E_p$ are effective (cf.\ \cite[III.6.1]{shaf}) and
$F_p$ is effective since $p_*|M_Z|=|M|$. 

For any rational number $c$ we can rearrange (\ref{noe.fano.ineq.thm}.3) to get
$$
\begin{array}{l}                                                                
K_Z+c|M_Z|\simq q^*\bigl(K_{X'}+c|M'|\bigr)+E_q \qtq{and}\\                                    
K_Z+c|M_Z|\simq p^*\bigl(K_{X}+c|M|\bigr)+E_p-cF_p.                                           
\end{array}                                                                     \eqno{(\ref{noe.fano.ineq.thm}.4)}
 $$
First we set 
$c=\tfrac1{m'}$.  Then $K_{X'}+\tfrac1{m'}|M'|\simq 0$, hence
$$                                                                               K_Z+\tfrac1{m'}|M_Z|\simq q^*\bigl(K_{X'}+\tfrac1{m'}|M'|\bigr)+E_q                           
\simq E_q\geq 0. 
%\eqno{(\ref{noe.fano.ineq.thm}.3)}
$$
Pushing this forward to $X$ we get that 
$$
K_{X}+\tfrac1{m'}|M|=p_*\bigl(K_Z+\tfrac1{m'}|M_Z|\bigr)\simq
p_*(E_q)\geq 0.
%\eqno{(\ref{noe.fano.ineq.thm}.4)}
$$
Since
$$
p_*(E_q)\simq K_{X}+\tfrac1{m'}|M|\simq K_X-\tfrac1{m'}mK_X=\tfrac{m-m'}{m'}(-K_X),
\eqno{(\ref{noe.fano.ineq.thm}.5)}
$$
we see that   $m \geq m'$.

Next set   $c=\tfrac1{m}$. Then we get that  
$$                                                                              
K_Z+\tfrac1{m}|M_Z|\simq p^*\bigl(K_{X}+\tfrac1{m}|M|\bigr)+E_p-\tfrac1{m}F_p                 
\simq E_p-\tfrac1{m}F_p.                                                        
%\eqno{(\ref{noe.fano.ineq.thm}.6)}
$$
Pushing this forward to $X$ yields
$$
K_{X'}+\tfrac1{m}|M'|= q_*\bigl(K_Z+\tfrac1{m}|M_Z|\bigr)\simq
q_*(E_p-\tfrac1{m}F_p).
%\eqno{(\ref{noe.fano.ineq.thm}.7)}
$$
As in (\ref{noe.fano.ineq.thm}.5) we obtain  that
$$
\tfrac{m'-m}{m}(-K_{X'})\simq K_{X'}+\tfrac1{m}|M'|\simq q_*(E_p-\tfrac1{m}F_p).
\eqno{(\ref{noe.fano.ineq.thm}.6)}
$$

{\it Basic alternative \ref{noe.fano.ineq.thm}.7.} 
\begin{itemize}
\item If $E_p-\tfrac1{m}F_p$ is not effective, then we declare  the linear system  $|M|$  to be  ``very singular.'' 
In our terminology,   $\bigl(X, \tfrac1m |M|\bigr)$  is not canonical. This is case (2).
\item  If $E_p-\tfrac1{m}F_p$ is  effective, then we declare  the linear system  $|M|$  to be  ``mildly singular.''  
In our terminology,   $\bigl(X, \tfrac1m |M|\bigr)$  is  canonical. We need to  prove that in this case $\Phi$ is an isomorphism.
\end{itemize}
\medskip

Thus assume from now on that   $E_p-\tfrac1{m}F_p$ is effective. Then  (\ref{noe.fano.ineq.thm}.6) implies that  $m'\geq m$.  Combining it with
(\ref{noe.fano.ineq.thm}.5) gives that $m'=m$ and then
(\ref{noe.fano.ineq.thm}.5) shows that $p_*(E_q)=0$. That is, $\supp E_q$ is $p$-exceptional. Since $X'$ is smooth, the support of $E_q$ is the whole
$q$-exceptional divisor $\ex(q)$. Thus every $q$-exceptional divisor is also
$p$-exceptional.

To see the converse, let $D\subset Z$ be an irreducible divisor that is not 
$q$-exceptional. Then $q_*(D)\simq r|M'|$ for some $r>0$. Thus
$$
r|M_Z|\simq q^*\bigl(r|M'|\bigr)\simq D+(\mbox{$q$-exceptional divisor}).
$$
Pushing forward to $X$ now gives that
$r|M|\simq p_*(D)$, since every $q$-exceptional divisor is also
$p$-exceptional. Here $p_*(D)\neq 0$ since $r>0$, so 
 $D$ is not  $p$-exceptional. This shows that
$\ex(p)=\ex(q)$. 

Finally set $Z:=p\bigl(\ex(p)\bigr)\subset X$, $Z':=q\bigl(\ex(q)\bigr)\subset X'$ and apply the following result of Matsusaka and Mumford \cite{mats-mumf}
to conclude that $\Phi$ is an isomorphism \qed

\begin{lem}\label{mat.mum.lem}
 Let $\Psi:Y\map Y'$ be a birational map between smooth projective varieties. Let $Z\subset Y$ and $Z'\subset Y'$ be closed sets of codimension $\geq 2$ such that $\Psi$ restricts to an isomorphism
$Y\setminus Z\cong Y'\setminus Z'$. Let $H$ be an ample divisor on $Y$ such that $H':=\Psi_*H$ is also ample.  
Then $\Psi$ is an isomorphism.
\end{lem}

Proof. We may assume that $H'$ and $H$ are both  very ample. Then
$$
|H'|=\bigl|H'_{Y'\setminus Z'}\bigr|=\Psi_*\bigl|H_{Y\setminus Z}\bigr|=\Psi_*|H|.
$$
Thus $\Psi_*|H|$ is base point free, hence  $\Psi^{-1}$ is everywhere defined. 
The same argument, with the roles of $Y,Y'$ reversed,  shows that $\Psi$ is also everywhere defined. So $\Psi$  is an isomorphism.\qed

\begin{rem} The proof of Theorem~\ref{noe.fano.ineq.thm}
also works if $X$ has canonical singularities, $X'$ has terminal singularities and they  both have  class number 1. 
\end{rem}

\section{Subvarieties of hypersurfaces}\label{sec.11}

Our aim is to prove that a subvariety of a smooth hypersurface can not be unexpectedly singular along a large dimensional subset. 
The claim and the method go back to  Fano and Segre; the first complete statement and proof is in \cite[Prop.5]{MR1970356}.

\begin{thm} \label{residual.mult.thm}
Let $X\subset \p^{n+1}$ be a smooth hypersurface,
$Z\subset X$ an irreducible subvariety 
  and 
$W\subset X$ a pure dimensional subscheme  such that
$\dim Z+\dim W\geq \dim X$. 
Assume that either $\dim Z<\dim W$ or $\dim Z=\dim W=\frac12 \dim X$ and $W$ is a complete intersection in $X$. 
Then
$$
\mult_ZW\leq \frac{\deg W}{\deg X}.
\eqno{(\ref{residual.mult.thm}.1)}
$$
\end{thm}

 We define the multiplicity $\mult_ZW$ in Paragraph~\ref{mult.say}. 
 See  \cite[Lem.13]{MR2171057} and  \cite[2.1]{MR3598821} for generalizations of the theorem to complete intersections.

\begin{say}[Proof of Step~\ref{main.steps}.2]
For part (a) set  $W:=D\in |mH|$. By B\'ezout's theorem, $\deg W=m\deg X$ so
$\mult_ZD\leq m$. For part (b) set  $W=M\cdot M$. Then $\deg W=m^2\deg X$ so
$\mult_Z(M\cdot M)\leq m^2$.  Note that $M\cdot M$ is a complete intersection in $X$, so the Theorem applies even if 
 $n=4$ and $\dim Z=2$.\qed
\end{say}

\begin{rem} 
The simplest special case of the theorem is when $W$ is an intersection of $X$ with a hyperplane. Then $\deg W=\deg X$ hence we claim that $W$   has only finitely many singular points. 
Equivalently,   a given hyperplane can be tangent to a smooth hypersurface only at finitely many points.
I encourage the reader to prove this; there are very easy  proofs but also messy ones. Note that this is truly a projective statement. For example, 
$(z-y^2x=0)$ is a smooth surface in $\a^3$ and the plane $z=0$ is tangent to it everywhere along the $x$-axis. 

Consider next the case when $W$ is an intersection of $X$ with a hypersurface of degree $d$. Then (\ref{residual.mult.thm}.1) says that
$W$ has multiplicity $\leq d$ at all but finitely many of its points ${\mathbf p}$. 
The easy geometric way to prove this would be to find a line $\ell$ in $X$ that passes through ${\mathbf p}$ but not contained in $W$. This sounds like  a reasonable plan if $\deg X\leq n$, since  in these cases there is a line through every point of $X$, see  \cite[V.4.3]{rc-book}, which also shows that  if $\deg X\geq 2n$ then a general $X$ does not contain any lines. 

In Proposition~\ref{residual.better.prop}, as replacements of lines,  we construct certain auxiliary subvarieties $Z^*$ that have surprisingly many intersections with $W$.

The extra assumption in case  $\dim Z=\dim W=\frac12 \dim X$ is necessary. Indeed, there are smooth hypersurfaces
$X\subset \p^{2n+1}$ that contain a linear space $L$ of dimension $n$.
Setting $Z=W=L$ we get that $\mult_ZW=1$ but $\frac{\deg W}{\deg X} =\frac{1}{\deg X}$. 
\end{rem}

\begin{say}[Multiplicity]\label{mult.say}  The simplest measure of a singularity is its {\it multiplicity.} Let $X=(h=0)\subset \a^n$ be an affine hypersurface and ${\mathbf p}=(p_1,\dots, p_n)$ a point on $X$. We can  write the  equation as
$$
h=\tsum  a_{i_1,\dots, i_n}(x_1-p_1)^{i_1}\cdots (x_n-p_n)^{i_n}.
$$
The multiplicity of $X$ at  ${\mathbf p}$,   denoted by $\mult_{\mathbf p}X$, is defined as
$$
\mult_{\mathbf p}X:=\min\{i_1+\cdots+ i_n: a_{i_1,\dots, i_n}\neq 0\}.
\eqno{(\ref{mult.say}.1)}
$$
The definition of multiplicity for other varieties is, unfortunately, more complicated.
Let $Y\subset  \a^n$ be a variety of dimension $m$  and ${\mathbf p}=(p_1,\dots, p_n)$ a point on $Y$. The following give the correct definition of the multiplicity $\mult_{\mathbf p}Y$, see \cite[Chap.5]{MR0453732} for details.
\begin{enumerate}\setcounter{enumi}{1}
\item Let $\pi:\a^n\to \a^{m+1}$ be a general projection. Then  $\pi(Y)$ is a hypersurface and $\mult_{\mathbf p}Y=\mult_{\pi(\mathbf p)}\pi(Y)$.
\item  If we are over $\c$, we can fix a small  Euclidean ball  $B(\epsilon)$ around ${\mathbf p}$, a general linear subspace $L$ of dimension $n-m$ through 
${\mathbf p}$ and count  the number of those intersection points 
of $X$ with a 
 general small translate of $L$ that are contained in $B(\epsilon)$.
\item The multiplicity also equals the limit
$$
\lim_{r\to \infty}  \tfrac{m!}{r^m} \dim_k k[x_1,\dots, x_n]/\bigl(I_Y, (x_1-p_1, \dots, x_n-p_n)^r\bigr).
$$
\end{enumerate}
 The first 2  are  old-style definitions that capture the essence but are not easy to work with rigorously, the 3rd is easy to use algebraically but
it is not even obvious that the limit exists; see \cite{MR0072531}.
The most complete modern treatment is given in \cite{fulton-int}.
Finally we set  
$$
\mult_ZY:=\min\{\mult_{\mathbf p}Y: {\mathbf p}\in Z\},
\eqno{(\ref{mult.say}.5)}
$$ 
and note that the minimum is achieved on a dense open subset.

We will also need the following.

\medskip
{\it Theorem \ref{mult.say}.6.}  Let $X\subset \p^{n+1}$ be a smooth hypersurface  and 
$Z, W\subset X$ irreducible subvarieties such that
$Z\cap W$ is finite and 
$\dim Z+\dim W=\dim X$. Assume furthermore that neither of them has dimension $\frac{n}{2}$.  Then
$$
\tsum_{\mathbf p} \mult_{\mathbf p}Z\cdot  \mult_{\mathbf p}W\leq 
\tfrac{\deg Z\cdot\deg W}{\deg X}.
$$

{\it Comments on the proof.} There are several theorems rolled into one here.

Intersection theory says that  if $X$ is any smooth projective variety   and 
$Z, W\subset X$ irreducible subvarieties such that $\dim Z+\dim W=\dim X$,
then they have a natural {\it intersection number,} denoted by   $(Z\cdot W)$.
Intersection theory can be developed completely algebraically, but working over $\c$ there is a shortcut. Both $Z,W$ have a homology class
$[Z]\in H_{2\dim Z}(X(\c), \z)$ and $[W]\in H_{2\dim W}(X(\c), \z)$ and then
$$
(Z\cdot W)=[Z]\cap [W]\in H_0(X(\c), \z)\cong \z.
\eqno{(\ref{mult.say}.7)}
$$ 
Furthermore, if $Z\cap W$ is finite then their   intersection number  $(Z\cdot W)$  is the sum of local terms, denoted by
 $(Z\cdot W)_{\mathbf p}$,  computed at each
${\mathbf p}\in Z\cap W$. Next  we need that
$$
(Z\cdot W)_{\mathbf p}\geq \mult_{\mathbf p}Z\cdot  \mult_{\mathbf p}W.
\eqno{(\ref{mult.say}.8)}
$$  
This very useful inequality does not seem to be included in  introductory
books.  It is  easy to derive it from \cite[Cor.A.14]{MR0453732}, 
see also
\cite[p.95]{MR0072531} or \cite[Cor.12.4]{fulton-int}.

Assume next that  $X\subset \p^{n+1}$ is a smooth hypersurface of degree $d$ and  $W$ is obtained as the intersection of $X$ by  $n-r$ hypersurfaces of degrees  $m_{r+1},\dots, m_n$.
If $\dim W=r$ then, by B\'ezout's theorem,  
$$
 \deg W=d\cdot m_{r+1}\cdots m_n 
\qtq{and}  (Z\cdot W)=\deg Z\cdot m_{r+1}\cdots m_n.
$$
Thus  we obtain that
 $$(Z\cdot W)=\tfrac{\deg Z\cdot\deg W}{\deg X}.
\eqno{(\ref{mult.say}.9)}
$$ 
It is not at all true that every $W$ can be obtained this way, but, by the 
  Lefschetz  hyperplane theorem (see \cite{MR0033557} or \cite[p.156]{gri-har}),  the  homology class of $W$ is  a rational multiple of a power of the hyperplane class, provided $\dim W\neq \tfrac{n}{2}$. Thus the above computation applies to every $W$ as in Theorem~\ref{residual.mult.thm}. 

 \end{say}

\begin{say}[Proof of Theorem~\ref{residual.mult.thm}]\label{residual.mult.thm.pf}
 If $\dim Z+\dim W> \dim X$  and the claim holds for all subvarieties
$Z'\subset Z$ of codimension 1 then it also holds for $Z$. Thus we may assume from now on that $\dim Z+\dim W= \dim X$.

Both the multiplicity and the degree is linear in irreducible components,
so write $W=\sum m_i W_i$ where the $W_i$ are irreducible subvarieties of $X$.

In Proposition~\ref{residual.better.prop} we construct a  subvariety $Z^*\subset X$  such that  $\dim Z^*=\dim Z$, 
  $\deg Z^*=(d-1)^r\deg Z$,
 $Z\cap Z^*$ consists of at least $(d-1)^r\deg Z$ distinct points  and
$W\cap Z^*$ is finite. 

% By assumption $\dim Z<\frac12 \dim X<\dim W$ thus
%  the intersection numbers $(W_i\cdot Z^*)$ are well defined.

There is nothing to prove if $Z\not\subset W_i$. 
Otherwise,
at each point of $Z\cap Z^*$ the intersection multiplicity of
$Z^*$ and $W_i$ is at least $\mult_ZW_i$  by (\ref{mult.say}.8). Therefore
$$
\bigl((d-1)^r\deg Z\bigr)\cdot \mult_ZW_i\leq  (W_i\cdot Z^*).
\eqno{(\ref{residual.mult.thm.pf}.1)}
$$
Next we use that 
$\tsum_i m_i(W_i\cdot Z^*)=(W\cdot Z^*)=\frac{\deg W\cdot \deg Z^*}{\deg X}$ by
(\ref{mult.say}.9).
Suming (\ref{residual.mult.thm.pf}.1) we get that 
$$
\bigl((d-1)^r\deg Z\bigr)\cdot \mult_ZW\leq (W\cdot Z^*)= \frac{\deg W\cdot (d-1)^r\deg Z}{\deg X}. 
\eqno{(\ref{residual.mult.thm.pf}.2)}
$$
Canceling $(d-1)^r\deg Z$ gives (\ref{residual.mult.thm}.1). \qed
\end{say}

Next we construct the subvariety $Z^*$ used in the above proof.

\begin{prop}\label{residual.better.prop} Let $X\subset \p^{n+1}$ be a smooth hypersurface of degree $d$. Let 
$Z\subset X$ be a subvariety of  dimension $r\leq \frac{n}{2}$ and $W_i\subset X$  a finite set of subvarieties.
Then there is a subvariety $Z^*$ of  dimension $r$ such that
\begin{enumerate}
\item  $\deg Z^*=(d-1)^r\deg Z$, 
\item $Z\cap Z^*$ consists of at least $(d-1)^r\deg Z$ distinct points, and
\item $\dim \bigl(Z^*\cap (W_i\setminus Z)\bigr)\leq \dim Z+\dim W_i-\dim X$ for every $i$.
\end{enumerate}
\end{prop}

The proof relies on the study of  certain  residual intersections.

\begin{say}[Residual intersection with cones] \label{res.int.defn}
Let $X\subset \p^{n+1}$ be a  hypersurface of degree $d$ and
$Z\subset X$ a  subvariety. Pick a point ${\mathbf v}\in \p^{n+1}$ and let $\langle{\mathbf v}, Z\rangle $ denote the cone over 
$Z$ with vertex ${\mathbf v}$, that is, the union of all lines
$\langle{\mathbf v}, z\rangle : z\in Z$.

If $\dim Z\leq n-1$ and ${\mathbf v}$ is general then $\langle{\mathbf v}, Z\rangle $
  has the same degree as $Z$ but 1 larger dimension.
If $\langle{\mathbf v}, Z\rangle $ is not contained in  $X$ then
$X\cap \langle{\mathbf v}, Z\rangle $ is a subscheme of $X$ of degree $=d\cdot \deg Z$. This subscheme contains $Z$, thus we can write
$$
X\cap \langle{\mathbf v}, Z\rangle=Z\cup Z^{\rm res}_{\mathbf v},
\eqno{(\ref{res.int.defn}.1)}
$$
where $Z^{\rm res}_{\mathbf v} $
 is called the {\it residual intersection} of the cone with $X$.  Note that
$$
\deg Z^{\rm res}_{\mathbf v}=(d-1)\cdot\deg Z.
\eqno{(\ref{res.int.defn}.2)}
$$
We are a little sloppy here; if $X$ is singular along $Z$ then
$Z^{\rm res}_{\mathbf v}$ is well defined as a cycle but  not  well defined as  a subscheme. 
We will always consider the case when $X$ is smooth at  general points  $z\in Z$ and ${\mathbf v}$ is not contained in the tangent plane of $X$ at 
$z$. If these hold then $\langle{\mathbf v}, Z\rangle $ is also smooth at
$z$ and hence $Z\not\subset Z^{\rm res}_{\mathbf v}$.
Our aim is to understand the intersection $Z\cap Z^{\rm res}_{\mathbf v}$.

Note that $Z\cap Z^{\rm res}_{\mathbf v}$ can be quite degenerate.
For example, let $X$ be the cone  $(x^n+y^n=z^n)\subset \p^3$ with vertex
at $(0{:}0{:}0{:}1)$ and 
$Z$ the line  $(x-z=y=0)$. Then $\langle{\mathbf v}, Z\rangle $ is a plane that contains $Z$, hence it contains the vertex of the cone. Thus
$X\cap \langle{\mathbf v}, Z\rangle$ is a union of $n$ lines through
 $(0{:}0{:}0{:}1)$. Thus $Z^{\rm res}_{\mathbf v}$ is a union of $n-1$ lines and
$Z\cap Z^{\rm res}_{\mathbf v}=(0{:}0{:}0{:}1)$, a single point. 

We see below that similar bad behavior does not happen for smooth hypersurfaces.
\end{say}

\begin{say}[Ramification linear system]\label{ram.lin.sys.defn}
 Let $X=(G=0)\subset \p^{n+1}$ be a hypersurface. The tangent plane  
$T_{\mathbf p}X$ at a smooth point
$(p_0{:}\cdots{:}p_{n+1})$ is given by the equation
$$
\tsum_i x_i \tfrac{\partial G}{\partial x_i}({\mathbf p})=0.
\eqno{(\ref{ram.lin.sys.defn}.1)}
$$
Let ${\mathbf v}:=(v_0{:}\cdots{:}v_{n+1})\in \p^{n+1}$ be a point  and
$\pi_{\mathbf v}:\p^{n+1}\map \p^n$ the projection from ${\mathbf v}$. 
The {\it ramification divisor} $R_{\mathbf v}$ of $\pi_{\mathbf v}|_X$ is the set of points whose tangent plane passes through ${\mathbf v}$. Thus
$$
R_{\mathbf v}=\bigl(\tsum_i v_i \tfrac{\partial G}{\partial x_i}=0\bigr)\cap X.
\eqno{(\ref{ram.lin.sys.defn}.2)}
$$
Thus the $|R_{\mathbf v}|$ form a linear system, called the
{\it ramification linear system,} 
which is the restriction of
the linear system of all first derivatives of  $G$. 
We denote it by $|R_X|$.  The base locus of  $|R_X|$ is exactly the singular locus $\sing X$.

Note that  $|R_X|\subset |(\deg X-1)H|_X$, where $H$ is the hyperplane class.

\end{say}

\begin{lem}\label{res.int.int.lem} Let $X\subset \p^{n+1}$ be a smooth hypersurface of degree $d$. Let 
$Z\subsetneq X$ be a subvariety of  dimension $r$ and $W_i\subset X$  a finite set of subvarieties.
Then, for general $ {\mathbf v}\in \p^{n+1}$,
\begin{enumerate}
\item  $Z^{\rm res}_{\mathbf v}\cap Z= R_{\mathbf v} \cap Z$ (set theoretically) and
\item  $Z^{\rm res}_{\mathbf v}\cap (W_i\setminus Z)$ has dimension
$\leq \dim Z+\dim W_i-n$.
\end{enumerate}
\end{lem}

Proof.  Set $\tau:=\pi_{\mathbf v}|_X$.
If  $\tau$ is unramified at ${\mathbf x}\in X$ then
it is a local isomorphism near ${\mathbf x}$, thus
$\langle{\mathbf v}, Z\rangle \cap X=\tau^{-1}\bigl(\tau(Z)\bigr)$
equals $Z$ near $X$. Thus 
$Z\cap Z^{\rm res}_{\mathbf v}\subset Z\cap R_{\mathbf v}$.
To see the converse, 
it is enough to prove that $Z\cap Z^{\rm res}_{\mathbf v}$ contains a dense open subset of $Z\cap R_{\mathbf v}$. Thus choose a point
 ${\mathbf x}\in Z$ that is smooth both on $X, Z$ and such that 
$\tau$ ramifies at ${\mathbf x}$ 
but $\tau|_Z$ does not. Then the vector pointing from
${\mathbf x}$ to ${\mathbf v}$ is also a tangent vector of
$\langle{\mathbf v}, Z\rangle \cap X$,  hence  ${\mathbf x}$ is a singular point of $\langle{\mathbf v}, Z\rangle \cap X$.  So
${\mathbf x}\in Z^{\rm res}_{\mathbf v}$, proving (1).

Note that  ${\mathbf p}\in Z^{\rm res}_{\mathbf v}\cap (W_i\setminus Z)$ iff a 
secant line connecting ${\mathbf p}$ with some  point of $Z$
passes through ${\mathbf v}$. 
The union of all secant lines connecting a point of $Z$ with a  different point of $W_i$ has dimension $\dim Z+\dim W_i+1$. Thus only a
$\dim Z+\dim W_i+1-(n+1)$ dimensional family of  secant lines passes through a general point of $\p^{n+1}$, proving (2).
\qed

\begin{say}[Proof of Proposition~\ref{residual.better.prop}]
\label{residual.better.prop.pf} 
Set $r=\dim Z$ and $Z_0:=Z$.  We inductively define
$$
Z_{i+1}:=(Z_i)^{\rm res}_{{\mathbf v}_i}\qtq{for general}  {\mathbf v}_i\in \p^{n+1}.
\eqno{(\ref{residual.better.prop.pf}.1)}
$$
We claim that $Z^*:=Z_r$ has the right properties. First note that
(\ref{residual.better.prop}.1) follows from (\ref{res.int.defn}.2).

Using (\ref{res.int.int.lem}.1)  $r$ times we see that
$Z\cap Z^*$ consists of the intersection points
$$
Z\cap  R_{{\mathbf v}_1}\cap \cdots\cap  R_{{\mathbf v}_r}
\eqno{(\ref{residual.better.prop.pf}.2)}
$$
for general ${\mathbf v}_i$. (If $r=\frac{n}{2}$, we may also get finitely many other points $Z_{i+1}\cap (Z\setminus Z_i)$; these we can ignore.)
 Since $X$ is smooth, $|R_X|$ is base point free, thus (\ref{residual.better.prop.pf}.2) consists of $(d-1)^r\deg Z$ points in general position.
(We use characteristic 0 at the last step.) \qed

\end{say}

\section{Multiplicity and canonical singularities}\label{sec.12}

One can usually  compute or at least estimate  the multiplicity of a divisor or a linear system at a point quite easily, thus it would be useful to be able decide using multiplicities whether  a pair $\bigl(X, c|M|\bigr)$
 is canonical or log canonical.
This turns out to be possible  for surfaces, less so for 3--folds, but the notions diverge more and more as the dimension grows.

If a pair $\bigl(X, c|M|\bigr)$ is not canonical, then there is a  non-canonical  exceptional divisor. 
We start with an example where this divisor is obtained by just one blow-up. Note that every  exceptional divisor can be obtained by repeatedly blowing up subvarieties, but the more blow-ups we need, the harder it is to connect the multiplicity with being canonical.

\begin{exmp}\label{mult.can.exmp.1}
 Let $X$ be a smooth variety, $Z\subset X$ a smooth subvariety of codimension $r$ and $|M|$ a linear system. Let $\pi:X'\to X$ denote the blow-up of $Z$ with exceptional divisor $E$. 
Then
$$
K_{X'}=\pi^*K_X+(r-1)E\qtq{and} \pi^*|M|=|M'|+\mult_Z|M|\cdot E.
$$
Thus 
$$
 K_{X'}+c|M'|\simq \pi^*\bigl(K_{X}+c|M|\bigr)+ \bigl(r-1-c\cdot \mult_Z|M|\bigr)E.
\eqno{(\ref{mult.can.exmp.1}.1)}
$$
Note that 
we can apply this to any subvariety, after we replace $X$ by $X\setminus \sing Z$. 
We have thus proved the following.

\medskip
{\it Claim \ref{mult.can.exmp.1}.2.} 
 Let $X$ be a smooth variety, $|M|$ a linear system and $Z\subset X$ a subvariety. Then the following hold.
\begin{enumerate}
\item[(a)] If $\bigl(X, c|M|\bigr)$ is canonical
then  $c\cdot \mult_Z|M|\leq \codim_XZ-1$.
\item[(b)] If $\bigl(X, c|M|\bigr)$ is log canonical
then  $c\cdot \mult_Z|M|\leq \codim_XZ$. \qed
\end{enumerate}
The problem we have is that the converse holds only for $n=2$ and only for part (a). Thus here our aim is to get some weaker converse statements in dimensions 2 and 3.
In order to do this, we need a good series of examples.

\medskip
{\it Claim \ref{mult.can.exmp.1}.3.} $\bigl(\a^n, c|\sum \lambda_i x_i^{m_i}|\bigr)$  is log canonical iff 
$$
c\leq \tfrac1{m_1}+\cdots + \tfrac1{m_n}. 
$$
A very useful way to think about this is the following.
If we assign weights to the variables  $w(x_i)=\tfrac1{m_i}$ then
the linear system becomes weighted homogeneous of weight 1. 
Thus, our condition says that
$$
c \cdot w\bigl(\tsum \lambda_i x_i^{m_i}\bigr)\leq w(x_1\cdots x_n).
\eqno{(\ref{mult.can.exmp.1}.4)}
$$
The claim is  easy to prove if all the $m_i$ are the same or if you know how to use weighted blow-ups, but can be very messy otherwise.
The case $n=2$ and $m_1=2$ is quite instructive and worth trying.

See \cite[Sec.6.5]{ksc} for details in general (using weighted blow-ups). 
\end{exmp}

The following lemma, which  
is  a partial converse to (\ref{mult.can.exmp.1}.2.a),  proves 
 Step~\ref{main.steps}.3.a.

\begin{lem}\label{main.steps.3.a.pf}
  Let $X$ be a smooth variety and $|M|$ a linear system. Assume that $c\cdot \mult_p|M|\leq 1$ for every point $p\in X$ and $\dim X\geq 2$. Then  
$\bigl(X, c|M|\bigr)$ is canonical.
\end{lem}

Proof. 
For one blow-up $\pi:X'\to X$ as in (\ref{mult.can.exmp.1}.1) we have the formula
$$
 K_{X'}+c|M'|\simq \pi^*\bigl(K_{X}+c|M|\bigr)+ \bigl(r-1-c\cdot \mult_Z|M|\bigr)E.
$$
Since $r\geq 2$, our assumption $c\cdot \mult_p|M|\leq 1$ implies that
$r-1-c\cdot \mult_Z|M|\geq 0$. 

If $\tau:X''\to X'$ is any birational morphism  and
$$
K_{X''}+c|M''|\simq \tau^*\bigl(K_{X'}+c|M'|\bigr)+ E'',
$$
then we get that
$$
K_{X''}+c|M''|\simq (\tau\circ\pi)^*\bigl(K_{X}+c|M|\bigr)+ E''+\bigl(r-1-c\cdot \mult_Z|M|\bigr)\tau^*E.
$$
If $\bigl(X', c|M'|\bigr)$ is canonical then $E''$ 
is  effective and so is $$E''+\bigl(r-1-c\cdot \mult_Z|M|\bigr)\tau^*E.$$
Thus  $\bigl(X, c|M|\bigr)$ is also canonical.
If $p'\in X'$ is any point and $p=\pi(p')$ then 
 $\mult_{p'}|M'|\leq \mult_p|M|$, thus 
$c\cdot \mult_{p'}|M'|\leq 1$ and we can use induction. 

The problem is that this seems to be an infinite induction, since we can keep blowing up forever. There are 2 ways of fixing this.

The easiest is to  use a log  resolution as in Definition~\ref{log.res.say} and stop when
the birational transform of $|M|$  becomes base point free, hence  canonical.

Theoretically it is better to focus on one divisor at a time and 
use a lemma of Zariski and Abhyankar, which is a very weak form of resolution;  
see  \cite[2.45]{km-book} or \cite[4.26]{ksc}. \qed

\medskip
{\it Remark \ref{main.steps.3.a.pf}.1.}  Another proof is the following. Let $p\in B\subset X$ be a general complete intersection curve. Then $c\cdot \bigl(|M|\cdot B\bigr)\leq 1$, hence
$\bigl(B, c |M|_B\bigr)$ is log canonical. By (\ref{connectedness.thm.cor})  this implies that 
$\bigl(X, c|M|\bigr)$ is canonical. \qed

\medskip

The following partial converse to 
(\ref{mult.can.exmp.1}.2.b)
is a reformulation of \cite{var76}, see also \cite[6.40]{ksc} for a proof.

\begin{thm}\label{lcthreshold.curve.formula}
Let $S$ be a  smooth surface  and $|M|$ a linear system such that
$p\in S$ is a non-log-canonical center of   $\bigl(S,c|M|\bigr)$.
Then one can choose local coordinates $(x,y)$ at $p$ and weights
$w(x)=a$ and  $w(y)=b$ such that
$$                                                                              |M|\subset\bigl|x^iy^j:  w(x^iy^j)> \tfrac1{c} w(xy)=\tfrac1{c}(a+b)\bigr|.                        
\qed
\eqno{(\ref{lcthreshold.curve.formula}.1)}
$$
\end{thm}

\begin{exmp} It can be quite hard to find the right coordinate system that works; it is frequently given by complicated power series. For example, \cite{MR533711} 
writes down  a  degree 6 polynomial
$g(x,y)$ that, in suitable   local coordinates becomes  $x^2+y^{20}$. 
(I do not doubt the claim but I have been unable to find a clear, non-computational  explanation.) Taking $a=10$ and $b=1$ shows that $\bigl(\a^2, c (g=0)\bigr)$ is
log canonical for $c\leq \tfrac{11}{20}$. 
Related bounds and examples are given in \cite{MR2795944}. 
\end{exmp}

The following consequence proves (\ref{main.steps}.3.2.b),
we derive  Step~\ref{main.steps}.3.b from it in Paragraph~\ref{end.pf.main.steps.3.b}.

\begin{cor}\cite{MR1798984}\label{lcthreshold.curve.formula.cor}
Let $S$ be a  smooth surface  and $|M|$ a movable linear system such that
$p\in S$ is a non-log-canonical point of   $\bigl(S,c|M|\bigr)$.
Then $(M\cdot M)_p>\frac{4}{c^2}$.
\end{cor}

{\it Remark.} Unlike for Lemma~\ref{main.steps.3.a.pf}, a direct induction does not seem to work, but \cite{MR1798984} sets up a more complicated inductive assumption and proves it one blow up at a time. The following argument, relying on Theorem~\ref{lcthreshold.curve.formula}, easily generalizes to all dimensions. (Unfortunately, this is less useful since 
Theorem~\ref{lcthreshold.curve.formula} does not generalize to higher dimensions.) 

\medskip
Proof. Assume first that in (\ref{lcthreshold.curve.formula}.1) we have $a=b$. 
Then every member of $|M|$ is a curve that has  multiplicity $>\frac{2}{c}$ at $p$ and
the intersection multiplicity is at least the product of the multiplicities.
(This is a special case of  (\ref{mult.say}.8), but it is much simpler; see 
\cite[IV.3.2]{shaf}.)
 Hence  the intersection multiplicity is $> \frac{4}{c^2}$.

In general we  get that members of $|M|$ have 
 multiplicity $
>\tfrac{1}{c}\bigl(1+\min\{\tfrac{a}{b},\tfrac{b}{a}\}\bigr)$ at $p$ 
and this only gives that $(M\cdot M)_p>\frac{1}{c^2}$.
Thus we need to equalize $a$ and $b$. The best way to do this is by a  weighted blow-up, see \cite[Sec.6.5]{ksc}, but here the following trick works.

After multiplying with the common denominator, we may assume that $a, b$ are integers.
Set $x=s^a$ and $y=t^b$. These define a degree $ab$ morphism
$\tau:\a^2_{st}\to \a^2_{xy}$.  The inclusion
$$                                                                              |M|\subset\bigl|x^iy^j:ai+bj>\tfrac1{c}(a+b)\bigr|                        
$$
of (\ref{lcthreshold.curve.formula}.1) is now transformed into
$$                                                                              \tau^*|M|\subset\bigl|s^{ai}t^{bj}:ai+bj>\tfrac1{c}(a+b)\bigr|  
 \subset\bigl|s^{m}t^{n}:m+n>\tfrac1{c}(a+b)\bigr|.                     
$$
That is, $\tau^*|M|$ has multiplicity $>\tfrac1{c}(a+b)$, hence
$(\tau^*M\cdot \tau^*M)_p>\tfrac{1}{c^2}(a+b)^2.$
Intersection multiplicities get multiplied by the degree of the map under pull-back, thus we conclude that
$(M\cdot M)_p>\tfrac{1}{c^2}\cdot \tfrac{(a+b)^2}{ab}\geq \tfrac{4}{c^2}$.\qed
\medskip

A 3-dimensional analog of Theorem~\ref{lcthreshold.curve.formula} 
was  conjectured in \cite{MR1798984}. The method of \cite{MR1311348} shows that it is a consequence of
a result of Kawakita \cite{MR1839287}. See also \cite[Chap.5]{ksc} for more details.

\begin{thm}\label{kawakita.thm}
Let $X$ be a  smooth threefold  and $|M|$ a linear system such that
$p\in X$ is a non-canonical center of   $\bigl(X,c|M|\bigr)$.
Then one can choose local coordinates $(x,y,z)$ at $p$ and weights
$w(x)=a$,  $w(y)=b$  and $w(z)=1$ such that
$$                                                                              |M|\subset\bigl|x^iy^jz^k:  w(x^iy^jz^k)>\tfrac1{c}w(xy)=\tfrac1{c}(a+b) \bigr|.           \qed             
$$
\end{thm}

\begin{say}[Summary]
Let $|M|:=|\sum \lambda_i g_i|$ be a linear system on $\a^n$.
\begin{itemize}
\item If $n=2$ then we can decide whether  $\bigl(\a^2, c|M|\bigr)$ is canonical at the origin just by looking at the degrees of the monomials that occur in the $g_i$.
\item If $n=2$ then we can decide whether  $\bigl(\a^2, c|M|\bigr)$ is log canonical  at the origin  by looking at the  monomials that occur in the $g_i$, provided we use the right coordinate system.
\item If $n=3$ then we can decide whether  $\bigl(\a^3, c|M|\bigr)$ is canonical  at the origin  by looking at the  monomials that occur in the $g_i$, provided we use the right coordinate system.
\item If $n\geq 4$ then the situation is more complicated, see Example~\ref{ksc.6.45}. However, as we discuss in Section~\ref{sec.nlc}, there is the following partial replacement.

We can frequently show that   $\bigl(\a^n, c|\sum \lambda_i g_i|\bigr)$ is not log canonical  at the origin  by looking at the monomials that occur in a Gr\"obner basis of  the ideal $(g_i)$.
\end{itemize}

\end{say}

\begin{exmp}\cite[6.45]{ksc}\label{ksc.6.45} For  $r\geq 5$ 
 consider the linear system
$$|M_r|:=\bigl|(x^2 +y^2 +z^2)^2, x^r , y^r, z^r\bigr|.
$$ Show that
$\bigl(\c^3, c|M_r|\bigr)$ is
 log canonical iff $c\leq \tfrac12+\tfrac1r$.
However, using coordinate  changes and weights only shows that $c\leq \frac34$.
% if $(x', y', z')$ is any coordinate system with weights $w$ then
% $w\bigl(|M_r|\bigr)\leq  \tfrac43 w(x' y' z')$.
\end{exmp}

\section{Hyperplane sections and canonical  singularities}\label{sec.13}

We start with the proof of Step~\ref{main.steps}.4.a.

\begin{say}[Bertini type theorems] \label{bertini.say}
The classical Berti theorem---for differentiable maps also known as Sard's theorem---says that a general member of a  base point free linear system  on a smooth variety is also smooth. This has numerous analogs, all saying that if a variety has certain types of singularities then a   general member of a  base point free linear system  also has only the same type of singularities.
Thus it is not surprising that the same holds for canonical and log canonical singularities. The log canonical case of the following  proves  Step~\ref{main.steps}.4.a.

\medskip
{\it Proposition \ref{bertini.say}.1.}  
Let $H\subset X$ be a  general member of a  base point free linear system  $|H|$. If  $\bigl(X, c|M|\bigr)$  is   canonical (resp.\ log canonical) then so is $\bigl(H, c|M|_H\bigr)$.
\medskip

Proof.  Choose a log resolution
$\pi:X'\to X$ as in  Definition~\ref{log.res.say} and write
$$
K_{X'}=\pi^*K_X+\tsum e_iE_i\qtq{and} \pi^*|M|=|M'|+\tsum a_iE_i,
\eqno{(\ref{bertini.say}.2)}
$$
where $|M'|$ is base point free and $\sum E_i$ has simple normal crossing singularities only. Thus
$$
K_{X'}+c|M'|\simq \pi^*\bigl(K_X+c|M|\bigr)+\tsum (e_i-ca_i)E_i.
\eqno{(\ref{bertini.say}.3)}
$$
  Note that $|H|$ gives us   base point free linear systems 
$|H'|$ on $X'$ and  $|H'|_{E_i}$ on  each $E_i$. 
The adjunction formula   (stated only for curves but proved in general in \cite[VI.1.4]{shaf})
says that
$K_{H}=\bigl(K_{X}+H\bigr)|_{H}$ and $K_{H'}=\bigl(K_{X'}+H'\bigr)|_{H'}$.
Adding $H'=\pi^*H$ to 
(\ref{bertini.say}.3) and restricting to $H$ and $H'$ 
we get  that
$$
K_{H'}+c|M'|_{H'}\simq \pi^*\bigl(K_H+c|M|_H\bigr)+\tsum (e_i-ca_i)(E_i\cap H'),
\eqno{(\ref{bertini.say}.4)}
$$
where $H'$ is  smooth and $\sum (E_i\cap H')$ has simple normal crossing singularities only. 
If $\bigl(X, c|M|\bigr)$ is canonical (resp.\ log canonical) then
$e_i-ca_i\geq 0$  (resp.\ $\geq -1$) for every $i$.
The  same $e_i-ca_i$ are involved in (\ref{bertini.say}.4), except that 
some of the $E_j\cap H'$ may be empty, in which case 
 $e_j-ca_j$  does not matter for 
 $\bigl(H, c|M|_H\bigr)$. In any case,  $\bigl(H, c|M|_H\bigr)$ is also canonical (resp.\ log canonical). \qed

\end{say}

Let us next see what happens if we try to use the same method to prove  Step~\ref{main.steps}.4.b.

\begin{say}\label{bertini.say.2}  Here we have a non-canonical center
$p\in X$ and 
we take an $H\in |H|$ that passes through the point $p$. If there is an exceptional divisor $E_j\subset X'$ such that $\pi(E_j)=\{p\}$, then
$\pi^*H\supset E_j$.  Hence $\pi^*H$ is not smooth, it is not even irreducible.
In this case we write
$$
\pi^*H=H'+\tsum m_i E_i.
\eqno{(\ref{bertini.say.2}.1)}
$$
Adding  $H'=\pi^*H-\sum m_i E_i$ to 
(\ref{bertini.say}.2) we get 
$$
K_{X'}+H'+c|M'|\simq \pi^*\bigl(K_X+H+c|M|\bigr)+\tsum (e_i-m_i-ca_i)E_i.
\eqno{(\ref{bertini.say.2}.2)}
$$
Thus  restricting (\ref{bertini.say.2}.2) to $H'$ and $H$ 
we  get  that
$$
K_{H'}+c|M'|_{H'}\simq \pi^*\bigl(K_H+c|M|_H\bigr)+\tsum (e_i-m_i-ca_i)(E_i\cap H').
\eqno{(\ref{bertini.say.2}.3)}
$$
At first sight we are done. If $p$ is a non-canonical center of
$\bigl(X, c|M|\bigr)$ then there is an $E_j$ such that
 $\pi(E_j)=\{p\}$ and $e_j-ca_j<0$. Since $H$ passes through $p$,
$m_j\geq 1$ also holds, so $ e_j-m_j-ca_j<-1$. Thus
$E_j\cap H'$ shows that $p$ is a non-log-canonical center of $\bigl(H, c|M|_H\bigr)$.

However, all this falls apart if $E_j\cap H'=\emptyset$. 
This can easily happen for some $E_j$, but it is enough to show that it can not happen for every $E_j$ for which $ e_j-m_j-ca_j<-1$.
This is what we discuss next.
\end{say}

The following 2 interconnected theorems have many names. In \cite{km-book} and \cite{kk-singbook} it is called {\it inversion of adjunction,} while \cite{puk-book}  uses 
{\it Shokurov-Koll\'ar connectedness principle.} A closely related result in complex analysis is the {\it Ohsawa-Takegoshi extension theorem,} proved in 
\cite{oh-ta}. 
The theorems were conjectured in \cite{sho-3ff} and proved in \cite[Sec.17]{k-etal}. 
The sharpest form was established in \cite{kawakita}, see also
\cite[Sec.4.1]{kk-singbook} for other generalizations.

For simplicity I state it only for smooth varieties, though the singular case is needed for most applications. The proof is actually a quite short application of Theorem~\ref{gen.kod.0.say}.3; see 
\cite[Sec.5.4]{km-book}  or  \cite[Chap.6]{ksc} for detailed treatments.

\begin{thm} \label{connectedness.thm} Let $X$ be a smooth variety and $\Delta$ an effective $\q$-divisor on $X$. Let $\pi:X'\to X$ be a proper, birational morphism and write
$$
K_{X'}\simq \pi^*(K_X+\Delta) +\tsum b_iB_i,
$$
where the $B_i$ are either $\pi$-exceptional or lie over $\supp \Delta$. 
Then   every fiber of
$$
\pi: \supp\bigl(\tsum_{b_i\leq -1} B_i\bigr)\to X
\qtq{is connected.} \qed
$$
\end{thm}

The following consequence is especially important.
The first part of it directly  implies Step~\ref{main.steps}.4.b--c,
the second part is also used in Section~\ref{sec.nlc}.

\begin{thm}\label{connectedness.thm.cor}
Let $X$ be a smooth variety and $|M|$ a linear system on $X$.
Let  $H\subset X$ be a smooth divisor.  Assume that
\begin{enumerate}
\item either $H$ contains  a non-canonical center $Z^{\rm c}$ of $\bigl(X, c|M|\bigr)$,
\item or $H$  has nonempty intersection with a non-log-canonical center $Z^{\rm lc}$ of $\bigl(X, c|M|\bigr)$.
\end{enumerate}
Then $\bigl(H, c|M|_H\bigr)$ is not log-canonical.
\end{thm}

Proof. Choose a log resolution $\pi:X'\to X$ and write
$\pi^*H=H'+\sum m_i E_i$. 
Choosing a general member $M\in |M|$ gives $M'\in |M'|$. 
We can  rearrange (\ref{bertini.say.2}.2) as
$$
K_{X'}\simq \pi^*\bigl(K_X+H+cM\bigr)-H'-cM'+\tsum (e_i-m_i-ca_i)E_i.
\eqno{(\ref{connectedness.thm.cor}.3)}
$$
Pick a point $p\in Z^{\rm c}$  (resp.\ $p\in H\cap Z^{\rm lc}$).
We may harmlessly assume that  $c<1$ (this is always the case for us) and  then   $-cM'$ does not contribute to the $\tsum_{b_i\leq -1} B_i$ in 
Theorem~\ref{connectedness.thm}. Thus we get  that
$$
F_p:=\pi^{-1}(p)\cap \bigl(H'\cup \tsum_{e_i-m_i-ca_i\leq -1} E_i\bigr)
\qtq{is connected.}
$$
If (1) holds then
there is an $E_j$ such that
 $p\in \pi(E_j)$,  $e_j-ca_j<0$ and $m_j\geq 1$.
If (2) holds then
there is an $E_j$ such that
 $p\in \pi(E_j)$ and  $e_j-ca_j<-1$. Thus, in both cases,
  $e_j-m_j-ca_j<-1$ and
$\pi^{-1}(p)\cap\tsum_{e_i-m_i-ca_i\leq -1} E_i$ is not empty. Since $F_p$ is connected,
we obtain that
$$
\pi^{-1}(p)\cap H'\cap \tsum_{e_i-m_i-ca_i\leq -1} E_i\neq\emptyset.
$$
Thus there is at least one divisor   $E_{j_0}$ such that 
 $$
 e_{j_0}-m_{j_0}-ca_{j_0}\leq -1\qtq{and}  p\in \pi\bigl(E_{j_0}\cap H'\bigr).
$$
Hence $E_{j_0}\cap H'$ gives the non-empty divisor that we needed in
(\ref{bertini.say.2}.3).
A small problem is that  we would like a strict inequality $e_{j_0}-m_{j_0}-ca_{j_0}< -1$. To achieve this, 
run the same argument with some $c'<c$. Then we get a $j_0$  such  that
$$
 e_{j_0}-m_{j_0}-ca_{j_0}<e_{j_0}-m_{j_0}-c'a_{j_0}\leq -1\qtq{and}  p\in \pi\bigl(E_{j_0}\cap H'\bigr).
$$
Thus  $\bigl(H, c|M|_H\bigr)$ is  not  log canonical.
This completes the proof of   Step~\ref{main.steps}.4.b. \qed

% \begin{say}[Proof of  Step~\ref{main.steps}.5]\label{main.steps.4.pf}
%  Let  $Y\subset \p^{n+1}$ be a smooth hypersurface,  $|M|\subset |mH|$  a movable linear system and $M\in |M|$ any  member. 
% Let $B_m\subset M$ be the set of points where $\mult_{\mathbf p}M> m$. This is a closed subset, thus it is either finite or it contains  a curve $C$. 
% The latter is impossible by    Theorem~\ref{residual.mult.thm}. 
% Thus $B_m$ is finite and   $\bigl(Y, \tfrac1m |M|\bigr)$ is   canonical outside $B_m$ by Lemma~\ref{main.steps.3.a.pf}. This proves Step~\ref{main.steps}.5.a.

% The argument for Step~\ref{main.steps}.5.b is similar. We use Theorem~\ref{residual.mult.thm} to obtain that there is a  1-dimensional closed subset $C_m$ such that $\mult_{\mathbf p}(M\cdot M)\leq m^2$ for
% every ${\mathbf p}\in Y\setminus C_m$. Then we claim that 
% $\bigl(Y, \tfrac1m |2M|\bigr)$ is  log canonical on $Y\setminus C_m$. 
% To see this  pick any ${\mathbf p}\in Y\setminus C_m$. The multiplicity is preserved by general hyperplane cuts, thus eventually we get a surface
% $S=H_1\cap\cdots\cap H_{n-2}$ containing ${\mathbf p}$ such that
% $\mult_{\mathbf p}(|M|_S\cdot |M|_S)\leq m^2$. Thus  $\bigl(S, \tfrac1m |2M|_S\bigr)$ is  log canonical at ${\mathbf p}$ by Corollary~\ref{lcthreshold.curve.formula.cor}. Finally we use
% Theorem~\ref{connectedness.thm.cor} to conclude that
%  $\bigl(Y, \tfrac1m |2M|\bigr)$ is also log canonical at ${\mathbf p}$.\qed

% \end{say} 

\begin{say}[End of the proof of  Step~\ref{main.steps}.3.b]
\label{end.pf.main.steps.3.b}
Let  $|M|$ be a movable linear system on a smooth variety $Y$ and
${\mathbf p}\in Y$ a non-log-canonical point of $\bigl(Y, c |M|\bigr)$. 
 The multiplicity at ${\mathbf p}$ is preserved by general hyperplane cuts through ${\mathbf p}$ 
and so is being non-log-canonical by Theorem~\ref{connectedness.thm.cor}.
Thus eventually we get a surface
$S=H_1\cap\cdots\cap H_{n-2}$ containing ${\mathbf p}$ 
such that  $\bigl(S, c |M|_S\bigr)$ is  non-log-canonical at ${\mathbf p}$.
Thus $\mult_{\mathbf p}(M\cdot M)=\mult_{\mathbf p}(|M|_S\cdot |M|_S)> 4/c^2$
by Corollary~\ref{lcthreshold.curve.formula.cor}.  \qed
\end{say}

\section{Global sections from isolated singularities}\label{sec.14}

%Following \cite{MR3455160, zhuang},
The proof of Step~\ref{main.steps}.6 is a combination of 4 lemmas, which are either quite easy to prove (Lemma~\ref{for.z.lem.1}) or have been well known
(Lemmas~\ref{mult.ideal.supp.lem}, \ref{for.z.lem.2} and \ref{for.z.lem.3}). Nonetheless, the power of their combination was not realized before \cite{zhuang}. 

We define the upper multiplier ideals ${\mathcal J}^+\bigl(c |N|\bigr)$ 
and ${\mathcal J}^+\bigl(I^c\bigr)$ 
for a
linear system $|N|$ and  an  ideal $I$ in Definition~\ref{mult.ideal.defn}.
We use the following of its properties.

\begin{lem}\label{mult.ideal.supp.lem}  Let $Y$ be a smooth variety and
$|N|$ a linear system. Then the support of $\o_Y/{\mathcal J}^+\bigl(c|N|\bigr) $ is the union of all non-log-canonical centers of
$\bigl(Y, c|N|\bigr)$.
\end{lem}

\begin{lem}\label{for.z.lem.1} Let $Y$ be a smooth variety and
$|N|$ a linear system.
Assume that $\bigl(Y, \tfrac{c}{2} |N|\bigr)$ is not log canonical. 
Then 
 $${\mathcal J}^+\bigl({\mathcal J}^+\bigl(c |N|\bigr)\bigr)\neq \o_Y.$$
\end{lem}

\begin{lem}\label{for.z.lem.2} Let $Y$ be a smooth variety and
$I\subset \o_Y$ an ideal sheaf that vanishes only at  finitely many points.
Assume that ${\mathcal J}^+(I)\neq \o_Y$. Then
$$
\dim \bigl(\o_Y/I\bigr)\geq \tfrac12 3^{\dim Y}.
$$
\end{lem}

\begin{lem}\label{for.z.lem.3} Let $Y$ be a smooth, projective variety, $H$ an ample divisor on $Y$ and $|N|$ a linear system such that $H\simq c|N|$.
Assume that $\bigl(Y, c |N|\bigr)$ is log canonical outside finitely many points. Then
$$
H^0\bigl(Y, \o_Y(K_Y+H)\bigr)\geq \dim \bigl(\o_Y/{\mathcal J}^+\bigl(c |N|\bigr)\bigr).
$$
\end{lem}

\begin{say}[Proof of Step~\ref{main.steps}.6 using Lemmas~\ref{mult.ideal.supp.lem}--\ref{for.z.lem.3}]

We apply Lemma~\ref{for.z.lem.1} to $|N|:=|2M|$ and $c=\frac1{m}$.
We get the ideal sheaf $I:={\mathcal J}^+\bigl(\tfrac1m |2M|\bigr)$
such that ${\mathcal J}^+(I)\neq \o_Y$. 
By Lemma~\ref{mult.ideal.supp.lem} $I$ vanishes only at  finitely many points.
Thus, by
 Lemma~\ref{for.z.lem.2},
$$
\dim \bigl(\o_Y/{\mathcal J}^+\bigl(\tfrac1m |2M|\bigr)\bigr)\geq \tfrac12 3^{\dim Y}.
$$
Finally using  Lemma~\ref{for.z.lem.3} for  $H:=2L$ says that
$$
H^0\bigl(Y, \o_Y(K_Y+2L)\bigr)\geq \dim \bigl(\o_Y/{\mathcal J}^+\bigl(\tfrac1m |2M|\bigr)\bigr)\geq \tfrac12 3^{\dim Y}.\qed
$$
\end{say}

We prove Lemma~\ref{mult.ideal.supp.lem} in Paragraph~\ref{mult.ideal.supp.lem.pf}. 
 Lemma~\ref{for.z.lem.2} is local at the points where $I$ vanishes, in fact, it is a quite general algebra statement about ideals. We discuss it in detail in Section \ref{sec.nlc}.
 Lemma~\ref{for.z.lem.3} is a restatement of Corollary~\ref{for.z.lem.3.rest};
 we explain its proof in  Section \ref{sec.vanish}.

\begin{say}[Proof of Lemma~\ref{for.z.lem.1}]\label{for.z.lem.1.pf}
Take a
log resolution $\pi:Y'\to Y$ as in  Definition~\ref{log.res.say}  and write
$$
K_{Y'}\sim \pi^*K_Y+\tsum e_i E_i\qtq{and} 
\pi^*|N|=|N'|+\tsum_i a_iE_i,
$$
where $|N'|$ is base point free.
By Definition~\ref{mult.ideal.say}, 
$$
{\mathcal J}^+\bigl(c|N|\bigr)=\pi_*\o_{Y'}\bigl(\tsum_i \rup{e_i-c'a_i}E_i\bigr),
$$
where $0<c-c'\ll 1$.  Thus if
$$
\pi^*{\mathcal J}^+\bigl(c|N|\bigr)=\o_{Y'}\bigl(-\tsum_i b_iE_i\bigr),
$$
then  $-b_i\leq \rup{e_i-c'a_i}$.
(We have to be a little careful here. We need to use a $\pi:Y'\to Y$
that is a log resolution for both  $|N|$ and ${\mathcal J}^+\bigl(c|N|\bigr)$.)
Therefore
$$
{\mathcal J}^+\bigl({\mathcal J}^+\bigl(c|N|\bigr)\bigr)
= \pi_*\o_{Y'}\bigl(\tsum_i \rup{e_i-(1-\epsilon)b_i}E_i\bigr),
\eqno{(\ref{for.z.lem.1.pf}.1)}
$$
for  $0<\epsilon\ll 1$.
If $b_i=0$ then $\rup{e_i-(1-\epsilon)b_i}=e_i-b_i$ and if
$b_i>0$ then 
$$\rup{e_i-(1-\epsilon)b_i}= e_i-b_i+\rup{\epsilon b_i}=e_i-b_i+1,
$$
since $b_i$ is an integer. Thus, in both cases
$$\rup{e_i-(1-\epsilon)b_i}\leq e_i-b_i+1\leq 2e_i+\rup{-c'a_i}+1\leq 
2e_i-ca_i+2.
$$
Since $\bigl(Y, \tfrac{c}{2} |N|\bigr)$ is not log canonical, 
there is an index $j$ such that
$e_j-\tfrac{c}{2}a_j<-1$. Then
$$
2e_j-ca_j+2=2\bigl(e_j-\tfrac{c}{2}a_j\bigr)+2<-2+2=0.
$$
Thus ${\mathcal J}^+\bigl({\mathcal J}^+\bigl(c|N|\bigr)\bigr)$ vanishes along
$\pi(E_j)$. \qed
\end{say}

\section{Review of vanishing theorems}\label{sec.vanish}

Here we prove  Lemma~\ref{for.z.lem.3}.
For this   we need to use the cohomology of coherent sheaves.
We use that the groups $H^i(Y, F)$ exist and that  a short exact sequence of sheaves leads to a long exact sequence of the cohomology groups.
For the uninitiated, \cite[Chap.B]{MR1442522} is a very good introduction.
 
We also need a vanishing theorem  which says that under certain assumptions the cohomology group $H^1(Y, F)$ is 0. The reader who is willing to believe Theorem~\ref{gen.kod.0.say.6} need not get into any further details. However, at first sight, the definition of the multiplier ideal  may appear rather strange, so I include an explanation of where these definitions and results  come from. 

\begin{say}[Vanishing and global sections]\label{van.glob.say}
Let $F$ be a coherent sheaf on a projective variety $Y$.
One way to estimate the dimension of $H^0(Y, F)$ from below is to 
identify a subsheaf  ${\mathcal S}(F)\into F$ and the corresponding quotient
$F\onto {\mathcal Q}(F)$,  write down the short exact sequence
$$
0\to  {\mathcal S}(F)\to F\to {\mathcal Q}(F)\to 0,
$$
and the beginning of its long exact sequence
$$
0\to H^0\bigl(Y,{\mathcal S}(F)\bigr)\to 
H^0\bigl(Y,F\bigr)\to H^0\bigl(Y,{\mathcal Q}(F)\bigr)\to 
H^1\bigl(Y,{\mathcal S}(F)\bigr).
$$
If the last term vanishes then
$$
\dim H^0\bigl(Y,F\bigr)\geq \dim  H^0\bigl(Y,{\mathcal Q}(F)\bigr).
$$
\end{say}

In our case we have a divisor $L$ and a linear system
$|M|$ such that  $c|M|\sim L$ for some $c$. We will use these data to construct
the subsheaf 
$$
{\mathcal S}\bigl(\o_Y(K_Y+L)\bigr)\subset \o_Y(K_Y+L),
$$
such that  a generalization of Kodaira's vanishing theorem applies to
${\mathcal S}\bigl(\o_Y(K_Y+L)\bigr) $.  
These vanishing theorems form a  powerful machine which  gives us a vanishing involving  $\o_Y(K_Y+L)$. (For our purposes in Section~\ref{sec.14},  a vanishing involving pretty much any other $\o_Y(aK_Y+bL)$ would be good enough, as long as $a, b$ are much smaller than $\dim Y$.) 

\begin{say}[Generalizations of Kodaira's vanishing theorem]\label{gen.kod.0.say}
 Kodaira's classical vanishing theorem  says that if $Y$ is a smooth, projective variety over $\c$ and $L$ an ample divisor on $X$ then
$H^i\bigl(Y, \o_Y(K_Y+L)\bigr)=0$ for $i>0$. 
It has various generalizations when $L$ is only close to being ample.

\medskip
{\it \ref{gen.kod.0.say}.1 (Close-to-ample divisors).} It turns out that Kodaira's vanishing theorem also works for a divisor $L$ if one can write it as
 $L\simq cA+\Delta$ where 
\begin{enumerate}
\item[(a)] $c>0$ and  $A$ is  nef and big (that is, $(A\cdot C)\geq 0$ for every curve $C\subset Y$  and $\bigl(A^{\dim Y}\bigr)>0$), and
\item[(b)]  $\Delta:=\sum d_iD_i$,  where $d_i\in [0,1)$ and $\sum D_i$ has 
simple  normal crossing singularities  only.
\end{enumerate}
In practice the condition that  $\sum D_i$ be a simple normal crossing divisor is very rarely satisfied, but log resolution (as in Definition~\ref{log.res.say}) allows us to reduce almost everything to this case. 
The basic vanishing theorem is the following, see \cite[Sec.2.5]{km-book} or 
\cite[9.1.18]{laz-book} for proofs.

\medskip
{\it Theorem \ref{gen.kod.0.say}.2 (Kawamata-Viehweg version).}  Let $X$ be a smooth, projective variety  and $L$ a  divisor as in (\ref{gen.kod.0.say}.1).   Then
$$
H^i\bigl(X, \o_X(K_X+L)\bigr)=0\qtq{for} i>0. \qed
$$  
\medskip

The following   versions are easy to derive from
(\ref{gen.kod.0.say}.1);  see \cite[2.68]{km-book}.

\medskip
{\it Theorem \ref{gen.kod.0.say}.3 (Grauert-Riemenschneider  version).}  Let $X$ be a smooth, projective variety, $\pi:X\to Y$ a birational morphism and $L$ a  divisor as in (\ref{gen.kod.0.say}.1). Then 
$$
R^i\pi_*\o_X(K_X+L)=0\qtq{for} i>0. \qed
$$  
\medskip

{\it Corollary \ref{gen.kod.0.say}.4.}  Let $X$ be a smooth, projective variety, $\pi:X\to Y$ a birational morphism and $L$ a  divisor  as in (\ref{gen.kod.0.say}.1). Then
$$
H^i\bigl(Y, \pi_*\o_X(K_X+L)\bigr)=0 \qtq{for} i>0.  \qed
$$
\end{say}

Next we show how we get vanishing theorems starting with
a linear system.

\begin{say}%[Multiplier ideal]
\label{mult.ideal.say}
 Let $Y$ be a smooth, projective variety over $\c$ and  $|M|$ a linear system on $Y$. Following (\ref{gen.kod.0.say}.1) we would like to get a 
nef and big divisor plus
  a divisor with simple normal crossing support.

Thus let $\pi:Y'\to Y$ be a log  resolution as in  Definition~\ref{log.res.say}.
 Write
$$
K_{Y'}=\pi^*K_Y+\tsum e_iE_i\qtq{and} \pi^*|M|=|M'|+\tsum a_i E_i.
\eqno{(\ref{mult.ideal.say}.1)}
$$
Thus the $E_i$ are either $\pi$-exceptional or belong to the base locus of
$\pi^*|M| $, and
we allow $e_i$ or $a_i$ to be 0.  Let $L$ be an ample divisor such that   $L\simq c|M|$. Then   we can write 
the pull-back of $K_Y+L $ as 
$$
\pi^*(K_Y+L)\simq K_{Y'}+c|M'|+\tsum (-e_i+ca_i)E_i.
\eqno{(\ref{mult.ideal.say}.2)}
$$
The right hand side starts to look like we could apply
(\ref{gen.kod.0.say}.2) to it, but there are 2 problems. 
 The coefficient
$(-e_i+ca_i)$ need not lie in the interval $[0,1)$ and, although  
$|M'|$ is nef, it need not be big. The latter can be arranged by keeping a little bit of $L$ unchanged. That is, pick 
$0<c'<c$  and  write 
the pull-back of $K_Y+L $ as 
$$
\pi^*(K_Y+L)\simq K_{Y'}+(c-c') \pi^*L+c'|M'|+\tsum (-e_i+c'a_i)E_i.
\eqno{(\ref{mult.ideal.say}.3)}
$$
Now  $\pi^*L+c'|M'|$ is nef and big, but the first problem remains.
Here we use that 
any number $a$ can be uniquely written as $a=\rdown{a}+\{a\}$ where
$\rdown{a}$ is an integer and $\{a\}\in [0,1)$. 
Furthermore, let  $\rup{a}:=-\rdown{-a}$ denote the rounding up of $a$. 
Applying this and
 rearranging we get that
$$
\begin{array}{l}
\pi^*(K_Y+L)+\tsum \rup{e_i-c'a_i}E_i\\
\qquad \simq K_{Y'}+(c-c') \pi^*L+c'|M'|+\tsum \{-e_i+c'a_i\}E_i.
\end{array}
\eqno{(\ref{mult.ideal.say}.4)}
$$
Now  the vanishing (\ref{gen.kod.0.say}.4) applies to the right hand side of 
(\ref{mult.ideal.say}.4).
Note also that
$$
\pi_*\o_{Y'}\bigl(\pi^*(K_X+L)+\tsum \rup{e_i-c'a_i}E_i\bigr)=
\o_Y(K_Y+L)\otimes \pi_*\o_{Y'}\bigl(\tsum \rup{e_i-c'a_i}E_i\bigr).
$$
This suggests that the basic object is  $\pi_*\o_{Y'}\bigl(\tsum \rup{e_i-c'a_i}E_i\bigr) $.
Note that this does not depend on $c'$  as long as $c-c'$ is small enough.
Indeed, then
$$
\rup{e_i-c'a_i} =
\left\{
\begin{array}{lcl}
\rup{e_i-ca_i} &\mbox{if}& e_i-ca_i\notin \z\qtq{and}\\
\rup{e_i-ca_i}+1 &\mbox{if}& e_i-ca_i\in \z.
\end{array}
\right.
$$
\end{say}

\begin{defn}[Multiplier ideal] \label{mult.ideal.defn}  Let $X$ be a smooth, projective variety over $\c$ and $|M|$ a linear system. The  {\it upper multiplier ideal} of
$c|M|$ is
$$
{\mathcal J}^+\bigl(c|M|\bigr):=\pi_* \o_{Y'}\bigl(\tsum \rup{e_i-c'a_i}E_i\bigr)
$$
for any $c'$ satisfying  $0<c-c'\ll 1$.
It is not hard to see that this does not depend on the choice of $\pi:Y'\to Y$.
Note that \cite[Sec.9.2]{laz-book} calls 
$${\mathcal J}\bigl(c|M|\bigr):=\pi_* \o_{Y'}\bigl(\tsum \rup{e_i-ca_i}E_i\bigr) $$ the multiplier ideal.  
Clearly 
${\mathcal J}^+\bigl(c|M|\bigr)={\mathcal J}\bigl(c'|M|\bigr)$
 for  $0<c-c'\ll 1$.

For an ideal sheaf $I(g_i)$ let  $|M|:=|\sum \lambda_i g_i|$ and set
${\mathcal J}^+\bigl(I^c\bigr):={\mathcal J}^+\bigl(c|M|\bigr)$. 
\end{defn}

\begin{say}[Proof of Lemma \ref{mult.ideal.supp.lem}]
\label{mult.ideal.supp.lem.pf} 
Set $W:=\supp \o_Y/{\mathcal J}^+\bigl(c|M|\bigr) $. Then $W$ is exactly the $\pi$-image of the support of the negative part of
$ \tsum \rup{e_i-c'a_i}E_i$.  If $e_i-ca_i<-1$ then 
$e_i-c'a_i<-1$ so $\rup{e_i-c'a_i}\leq -1$. 
If $e_i-ca_i\geq -1$  and $a_i>0$ then 
$e_i-c'a_i>-1$ so $\rup{e_i-c'a_i}\geq 0$.
If $a_i= 0$ then $e_i-ca_i=e_i\geq 0$. \qed
\end{say}

Now we can apply (\ref{gen.kod.0.say}.4) and get the following, see 
\cite[Sec.9.4]{laz-book}. 

\begin{thm}[Nadel vanishing]\label{gen.kod.0.say.6} Let $Y$ be a smooth, projective variety, $L$ an ample divisor and  $|M|$ a linear system 
such that $c|M|\simq L$. Then 
$$
H^i\bigl(Y, \o_Y(K_Y+L)\otimes {\mathcal J}^+\bigl(c|M|\bigr)\bigr)=0
\qtq{for } i>0.  \qed
$$ 
\end{thm}

As we discussed in Paragraph~\ref{van.glob.say}, this immediately implies the following.

\begin{cor}\label{for.z.lem.3.rest} Let $Y$ be a smooth, projective variety, $L$ an ample divisor on $Y$ and $|M|$ a linear system such that $L\simq c|M|$.
Assume that $\bigl(Y, c |M|\bigr)$ is log canonical outside finitely many points. Then
$$
H^0\bigl(Y, \o_Y(K_Y+L)\bigr)\geq \dim \bigl(\o_Y/{\mathcal J}^+\bigl(c |M|\bigr)\bigr). \qed
$$
\end{cor}

\section{Review of monomial ideals}\label{sec.nlc}

In this section we prove  Lemma~\ref{for.z.lem.2}. Its claim  is local at the points where $I$ vanishes, we can thus work using  local coordinates at a point. Though not important, it is notationally simpler to pretend that we work at the origin of $\a^n$.  (This is in fact completely correct, one needs to argue that $Y$ and $\a^n$ have the same completions, \cite[Sec.II.2.2]{shaf}.)

As a general rule, an ideal is log canonical iff it contains  low multiplicity polynomials. In this section we give a precise version of this claim. Key special cases of the following are proved  by
Reid \cite{r-c3f} and Corti \cite{MR1798984}. More general
versions are in \cite{MR2047683, MR1828466}. An excellent  detailed treatment of this topic is given in \cite[Chap.9]{laz-book}, so I concentrate on the  the definitions and explanations, leaving the details to \cite{laz-book}.

\begin{thm}\label{col.nonlc.ideal.thm}
Let $I\subset R:=k[x_1,\dots, x_n]$ be an ideal vanishing only at the origin.
Assume that $I$  is not log canonical. Then
$$
\dim \bigl(R/I)\geq \min_{a_1,\dots, a_n\geq 0} \#\Bigl\{ {\mathbb N}^n\cap \bigl(
\tsum a_ir_i\leq \tsum a_i\bigr)\Bigr\}.
$$
\end{thm}

The proof is given in 2 steps. We first reduce to the case of monomial ideals in (\ref{def.to.mon.say}.4) and then  to counting lattice points in a simplex (\ref{mon.ideal.say}.2). 
Following the proof shows that the lower bound is sharp, but I do not know a closed formula for it. However,
a simple argument, given in Paragraph~\ref{latt.pt.count},  gives the following.

\begin{cor}\label{col.nonlc.ideal.cor}
Let $I\subset k[x_1,\dots, x_n]$ be a non-log-canonical ideal that vanishes only at the origin. Then
$$
\dim \bigl(R/I)\geq \tfrac12 3^n.
$$
\end{cor}

\begin{say}[Deformation to monomial ideals]\label{def.to.mon.say}
(See \cite[Chap.2]{clo-book} for details.)
Let $I\subset R:=k[x_1,\dots, x_n]$ be an ideal.
Write every $g\in R$ as  
$g=\operatorname{in}(g)+\operatorname{rem}(g)$ where
$\operatorname{in}(g):=a_g \tprod x_i^{r_i} $
 is the lexicographically lowest monomial that appears in $g$ with nonzero coefficient.  Define the {\it initial ideal} of $I$  (with respect to the lexicographic ordering) as
$$
\operatorname{in}(I):=\bigl( \operatorname{in}(g): g\in I\bigr).
\eqno{(\ref{def.to.mon.say}.1)}
$$
Thus $ \operatorname{in}(I)$ is generated by monomials and
it is not hard to see that
$$
\dim (R/I)=\dim \bigl(R/\operatorname{in}(I)\bigr).
\eqno{(\ref{def.to.mon.say}.2)}
$$
A key property is the following.

\medskip
{\it Proposition \ref{def.to.mon.say}.3.} If $\operatorname{in}(I)$ is  log canonical then so is  $I$. 
\medskip

{\it Comments on the proof.} Choose  integers $1\leq w_1\ll \cdots\ll w_n$.
For $g\in R$ let $w(g)$ denote the largest $t$ power that divides
$g(t^{w_1}x_1,\dots,  t^{w_n}x_n)$. Then
$$
t^{-w(g)}g(t^{w_1}x_1,\dots,  t^{w_n}x_n)=
\operatorname{in}(g)(x_1,\dots,  x_n)+t(\mbox{other terms}).
\eqno{(\ref{def.to.mon.say}.3.a)}
$$
Any finite collection of these defines a 
linear system  $|M|$ on $Y:=\a^{n+1}$ with coordinates
$(x_1,\dots, x_n,t)$.

If we choose $w_i$ that work for a Gr\"obner basis  $g_i\in I$, then we get $|M|$  whose restriction  to $(t=0)$ gives
$I_0=\operatorname{in}(I)$ and  to $(t=\lambda)$  gives 
$I_{\lambda}\cong I$ for $\lambda\neq 0$. 

If $I_0$ is log canonical then 
 $\bigl(Y, |M|\bigr)$ is also  log canonical by  Theorem~\ref{connectedness.thm.cor}.2, and so is 
 $I_{\lambda}\cong I$  by  (\ref{bertini.say}.1).  \qed

\medskip

Combining (\ref{def.to.mon.say}.2) and (\ref{def.to.mon.say}.3)
gives the first reduction step of the proof of Theorem~\ref{col.nonlc.ideal.thm}. 

\medskip
{\it Corollary \ref{def.to.mon.say}.4.}
If Theorem~\ref{col.nonlc.ideal.thm} holds for monomial ideals then it also holds for all ideals. \qed
\end{say}

\begin{say}[Monomial ideals]\label{mon.ideal.say}
Let $I\subset k[_1,\dots, x_n]$ be a {\it monomial ideal,} that is, an ideal generated by monomials. A very good description of $I$ is given by its {\it Newton polytope.} 

For $\prod x_i^{r_i}\in I$ we mark   the point $(r_1,\dots, r_n)$ with a big dot for elements of $I\setminus (x_1,\dots, x_n)I$  and with an invisible dot for elements of $(x_1,\dots, x_n)I$.
The   Newton polytope
is the boundary of the convex hull of the
marked  points, as in the next example.
$$
\begin{picture}(100,80)(0,20)

\put(0,0){\line(1,0){80}}
\put(0,0){\line(0,1){80}}

\put(0,0){\circle*{2}}
\put(0,10){\circle*{2}}
\put(0,20){\circle*{2}}
\put(0,30){\circle*{2}}
\put(0,40){\circle*{2}}
\put(0,50){\circle*{2}}
\put(0,60){\circle*{2}}
\put(0,70){\circle*{2}}

\put(10,0){\circle*{2}}
\put(20,0){\circle*{2}}
\put(30,0){\circle*{2}}
\put(40,0){\circle*{2}}
\put(50,0){\circle*{2}}
\put(60,0){\circle*{2}}
\put(70,0){\circle*{2}}
\put(80,0){\circle*{2}}

\put(0,70){\circle*{4}}
\put(10,50){\circle*{4}}
\put(20,30){\circle*{4}}
\put(40,10){\circle*{4}}
%\put(50,30){\circle*{4}}
\put(60,0){\circle*{4}}

\put(0,70){\line(1,-2){20}}
\put(20,30){\line(1,-1){20}}
\put(40,10){\line(2,-1){20}}
\end{picture}
\begin{array}{c}
\mbox{The Newton polygon of}\\
 (y^7, y^5x, y^3x^2, yx^4,x^6).\\
{\ }
\end{array}
$$
\medskip

A face of the Newton polytope is called {\it central}  if it contains a point all of whose coordinates are equal. 
$$
\begin{picture}(100,80)(0,20)
\put(0,0){\line(1,0){80}}
\put(0,0){\line(0,1){80}}

\put(0,0){\circle*{2}}
\put(0,10){\circle*{2}}
\put(0,20){\circle*{2}}
\put(0,30){\circle*{2}}
\put(0,40){\circle*{2}}
\put(0,50){\circle*{2}}
\put(0,60){\circle*{2}}
\put(0,70){\circle*{2}}

\put(10,0){\circle*{2}}
\put(20,0){\circle*{2}}
\put(30,0){\circle*{2}}
\put(40,0){\circle*{2}}
\put(50,0){\circle*{2}}
\put(60,0){\circle*{2}}
\put(70,0){\circle*{2}}
\put(80,0){\circle*{2}}

\put(0,70){\circle*{4}}
\put(10,50){\circle*{4}}
\put(20,30){\circle*{4}}
\put(40,10){\circle*{4}}
%\put(50,30){\circle*{4}}
\put(60,0){\circle*{4}}

\put(0,70){\line(1,-2){20}}
\put(0,50){\line(1,-1){50}}
\put(40,10){\line(2,-1){20}}
\end{picture}
\begin{array}{c}
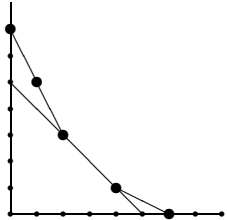
\mbox{The Newton polygon of}\\
 (y^7, y^5x, y^3x^2, yx^4,x^6),\\
\mbox{with central face extended.}\\
{\ }
\end{array}
$$

\medskip

The next version of (\ref{mult.can.exmp.1}.3)  follows from \cite{r-c3f}; see \cite{MR1828466} for various generalizations and \cite[Sec.9.3.C]{laz-book} for proof.

\medskip
{\it Proposition \ref{mon.ideal.say}.1.} A monomial ideal $I$ is log canonical iff its Newton polytope contains the point $(1,\dots, 1)$.  \qed
\medskip

Thus  $I$ is not log canonical  iff a central face of  its Newton polytope contains a point $(d,\dots, d)$ with $d>1$. The equation of this face can then be written as
$$
\tsum a_i r_i= d\tsum a_i\qtq{for some}  a_i\geq  0.
$$
In particular, a monomial $\prod x_i^{r_i}$ is not contained in $I$ if
$\tsum a_i r_i\leq \tsum a_i$. We have thus proved the following.

\medskip
{\it Corollary \ref{mon.ideal.say}.2.} A monomial ideal $I$ is not log canonical iff 
there is a simplex
$$\Delta({\mathbf a}):=\bigl(0\leq r_i,\ \tsum a_i r_i\leq \tsum a_i\bigr)$$
that is disjoint from the Newton polytope of $I$.
If this holds then
$$
\dim \bigl(R/I\bigr)\geq 
\left(
\begin{array}{c}
\mbox{number of lattice points}\\
\mbox{in the simplex $\Delta({\mathbf a})$}
\end{array}
\right). \qed
$$
\end{say}

\begin{say}[Lattice points in simplices]\label{latt.pt.count}
We thus need to estimate from below the number of  lattice points in the $n$-simplex
$\bigl(0\leq r_i, \tsum a_i r_i\leq \tsum a_i\bigr)$, independent of the $a_i$.
I could not find the optimal values.

The lower bound  $\tfrac12 3^n$ comes from the observation that if
$r_i\in \{0,1,2\}$ then either
$(r_1,\dots, r_n)$  or $(2-r_1,\dots, 2-r_n)$ satisfies  $\tsum a_i r_i\leq \tsum a_i$.   We can do  a little better by adding the points with coordinates $3,\dots, n$ on at least 1 of the coordinate axes.

Another  lower bound is
 $\tfrac1{n}\tbinom{2n}{n}$, which is asymptotically  $4^n/(n\sqrt{\pi n})$.
 This comes from the observation  that  if
$\sum r_i\leq n$ then at least one of the cyclic permutations of
$(r_1,\dots, r_n)$ satisfies  $\tsum a_i r_i\leq \tsum a_i$.

The first bound  is better for  $n\leq 5$, the second for $n\geq 6$. 

We can also combine the 2 bounds to get
$$
\tfrac1{n}\bigl[\tbinom{2n}{n}-\tfrac12 3^n\bigr]+\tfrac12 3^n=
\tfrac1{n}\tbinom{2n}{n}+\tfrac{n-1}{2n}3^n.
$$
\end{say}

%\bibliography{refs}
\def\cprime{$'$} \def\cprime{$'$} \def\cprime{$'$} \def\cprime{$'$}
  \def\cprime{$'$} \def\cprime{$'$} \def\dbar{\leavevmode\hbox to
  0pt{\hskip.2ex \accent"16\hss}d} \def\cprime{$'$} \def\cprime{$'$}
  \def\polhk#1{\setbox0=\hbox{#1}{\ooalign{\hidewidth
  \lower1.5ex\hbox{`}\hidewidth\crcr\unhbox0}}} \def\cprime{$'$}
  \def\cprime{$'$} \def\cprime{$'$} \def\cprime{$'$}
  \def\polhk#1{\setbox0=\hbox{#1}{\ooalign{\hidewidth
  \lower1.5ex\hbox{`}\hidewidth\crcr\unhbox0}}} \def\cdprime{$''$}
  \def\cprime{$'$} \def\cprime{$'$} \def\cprime{$'$} \def\cprime{$'$}
\providecommand{\bysame}{\leavevmode\hbox to3em{\hrulefill}\thinspace}
\providecommand{\MR}{\relax\ifhmode\unskip\space\fi MR }
\providecommand{\MRhref}[2]{%
  \href{http://www.ams.org/mathscinet-getitem?mr=#1}{#2}
}
\providecommand{\href}[2]{#2}

%\bigskip

\noindent  Princeton University, Princeton NJ 08544-1000

{\begin{verbatim} kollar@math.princeton.edu\end{verbatim}}

\end{document}